\documentclass[11pt]{amsart}

\usepackage{epsfig, subfig, amsmath, amssymb, latexsym, amsfonts, xcolor, easybmat, bbm}

\newcommand\ol{\ensuremath{\overline}}

\newcommand\eop{{{\hfil \ensuremath \Box}}}

\newcommand\eps{\ensuremath {\varepsilon}}

\newcommand\dss{\ensuremath{\displaystyle}}

\newenvironment{cor}{\subsection{}{\textbf {Corollary.}}\em}{}
\newenvironment{defn}{\subsection{}{\textbf {Definition.}}\em}{\smallskip}
\newenvironment{eg}{\subsection{}{\textbf {Example.}}}{\smallskip}

\newenvironment{lem}{\subsection{}{\textbf {Lemma.}}\em}{}
\newenvironment{notation}{\subsection{}{\textbf{Notation.}}}{\smallskip}

\newenvironment{prop}{\subsection{}{\textbf {Proposition.}}\em}{}
\newenvironment{rems}{\subsection{}{\textbf {Remarks.}}}{\smallskip}
\newenvironment{rem}{\subsection{}{\textbf {Remark.}}}{\smallskip}
\newenvironment{thm}{\subsection{}{\textbf {Theorem.}}\em}{}

\newenvironment{pf}{\noindent{\textbf {Proof.}}} {\begin{flushright}\eop \end{flushright}\smallskip}

\newcommand\cB{\ensuremath{\mathcal B}}

\newcommand\cD{\ensuremath{\mathcal D}}
\newcommand\cE{\ensuremath{\mathcal E}}

\newcommand\cH{\ensuremath{\mathcal H}}

\newcommand\cK{\ensuremath{\mathcal K}}

\newcommand\cM{\ensuremath{\mathcal M}}
\newcommand\cN{\ensuremath{\mathcal N}}

\newcommand\cS{\ensuremath{\mathcal S}}

\newcommand\cU{\ensuremath{\mathcal U}}
\newcommand\cV{\ensuremath{\mathcal V}}
\newcommand\cW{\ensuremath{\mathcal W}}
\newcommand\cX{\ensuremath{\mathcal X}}

\newcommand\bbC{\ensuremath{\mathbb C}}

\newcommand\bbM{\ensuremath{\mathbb M}}
\newcommand\bbN{\ensuremath{\mathbb N}}

\newcommand\bbR{\ensuremath{\mathbb R}}

\newcommand\bbT{\ensuremath{\mathbb T}}

\newcommand\bbZ{\ensuremath{\mathbb Z}}

\newcommand\ttt{\ensuremath{\textsc}}

\newcommand\bofh{\ensuremath{\cB ( \cH)}}

\newcommand\hilb{\ensuremath{\mathcal H}}

\newcommand\norm{\ensuremath {\Vert}}

\newcommand\rank{\ensuremath{\mathrm{rank}\, }}

\newcounter{asst}
\newcounter{lab}
\newcounter{asstAA}
\newcounter{asstBB}
\newcounter{asstCC}
\newcounter{asstDD}
\newcounter{asstEE}
\newcounter{asstFF}
\newcounter{asstGG}
\newcounter{asstHH}
\newcounter{asstII}
\newcounter{asstJJ}
\newcounter{asstKK}
\newcounter{asstLL}
\newcounter{asstMM}
\newcounter{suppAA}
\definecolor{myred}{rgb}{0.6,0,0}
\definecolor{myblue}{rgb}{0,0.2,0.4}
\definecolor{mygreen}{rgb}{0.2,0.6, 0.5}

\newcommand\symone{\ensuremath{\ttt{Sym}_1}}
\newcommand\symtwo{\ensuremath{\ttt{Sym}_2}}
\newcommand\symthree{\ensuremath{\ttt{Sym}_3}}
\newcommand\symfour{\ensuremath{\ttt{Sym}_4}}
\newcommand\symk{\ensuremath{\ttt{Sym}_k}}

\newcommand\uncom{\ensuremath{\Gamma_\cU}}

\newcommand\tr{\ensuremath{\textsc{tr}\,}}

\begin{document}


\title{On products of symmetries acting on Hilbert spaces}


\thanks{${}^1$ Research supported in part by NSERC (Canada).}
\thanks{${}^2$ Research supported in part by National Natural Science Foundation of China (No.: 12471123).}

\thanks{{\ifcase\month\or Jan.\or Feb.\or March\or April\or May\or
June\or
July\or Aug.\or Sept.\or Oct.\or Nov.\or Dec.\fi\space \number\day,
\number\year}}
\author
	[L.W. Marcoux]{{Laurent W.~Marcoux${}^1$}}
\address
	{Department of Pure Mathematics\\
	University of Waterloo\\
	Waterloo, Ontario \\
	Canada  \ \ \ N2L 3G1}
\email{LWMarcoux@math.uwaterloo.ca}


\author
	[H. Radjavi]{{Heydar Radjavi}}
\address
	{Department of Pure Mathematics\\
	University of Waterloo\\
	Waterloo, Ontario \\
	Canada  \ \ \ N2L 3G1}
\email{hradjavi@uwaterloo.ca}


\author
	[Y.H.~Zhang]{{Yuanhang~Zhang${}^2$}}
\address
	{School of Mathematics\\
	Jilin University\\
	Changchun 130012\\
	P.R. CHINA}
\email{zhangyuanhang@jlu.edu.cn}


\begin{abstract}
Let $\hilb$ be a complex, separable Hilbert space (of finite or infinite dimension), and let $\cU(\hilb)$ denote the group of unitary operators on $\hilb$.   A symmetry is, by definition, a unitary operator $J$ with $J^2 =I$.  Denote by $\symk(\hilb)$ the subset of $\cU(\hilb)$ consisting of those operators expressible as a product of $k$ symmetries.   It is known that $\cU(\hilb) = \symfour(\hilb)$ if $\dim \, \hilb = \infty$, while the only additional condition in finite dimensions is that the determinant be $\pm 1$.   Of all the sets $\symk(\hilb)$ with $k \in \{ 1, 2, 3, 4\}$, the case $k =3$ has been the most stubborn to characterise.   Among other things, we investigate which elements of $\symthree(\hilb)$ possess exactly two eigenvalues in the setting where $\hilb$ is finite-dimensional.  We also consider the problem:  when is the unitary orbit of an operator $T$, i.e., the set
\[
\{ U^* T U : U \in \cU(\hilb) \} \]
the same as its $\symk$-orbit, i.e., the set
\[
\{ U^* T U: U \in \symk(\hilb)\}?\]
Clearly, the cases of interest are when $k \le 3$.
\end{abstract}





\maketitle
\markboth{\textsc{  }}{\textsc{}}


\section{Introduction} \label{section01}


\subsection{} \label{sec1.01}
Let $\hilb$ be a complex Hilbert space, and denote by $\bofh$   the algebra of bounded linear operators acting on $\hilb$.  We shall write $\cU(\hilb) = \{ U \in \bofh: U \text{ is unitary}\}$ for the set of all unitary operators acting on $\hilb$.  For $\cX  \subseteq \bofh$, $\cU(\cX) = \{ U^* X U: X \in \cX, U \in \cU(\hilb)\}$ denotes the unitary orbit of $\cX$.  When $\cX = \{ X \}$ is a singleton set, we abbreviate the notation to $\cU(X)$.  If $\cU(\cX) \subseteq \cX$, we say that $\cX$ is \textbf{invariant under unitary conjugation}.
An operator $J \in \bofh$ is said to be a \textbf{symmetry} if $J = J^* = J^{-1}$.   It was R.V.~Kadison who first raised the question of which operators may be expressed as a finite product of symmetries~\cite{HalmosKakutani1958}.  Since each symmetry is unitary, any such product  must result in a unitary operator.


\begin{notation} \label{not1.02}
For $k \in \bbN$, set
\[
\symk(\hilb) = \{ V \in \cB(\hilb) : V = J_1 J_2 \cdots J_k : J_j \text{ is a symmetry}, 1 \le j \le k \}. \]
\end{notation}

The set  $\symone(\hilb)$ of symmetries on $\hilb$  is easily seen to be invariant under unitary conjugation, from which it in turn  follows that $V \in \symk(\hilb)$ if and only if $\cU(V) \subseteq \symk(\hilb)$, $k \ge 1$.  A classical result of Halmos and Kakutani shows that we need only consider $\symk(\hilb)$ for $1 \le k \le 4$.


\begin{thm} \cite{HalmosKakutani1958}  \label{HalmosKakutani1958}
If $\hilb$ is an infinite-dimensional Hilbert space, then every unitary operator on $\hilb$ is a product of four symmetries.   That is,
\[
\cU(\hilb) = \symfour(\hilb). \]
Furthermore, if $\omega$ is a primitive third root of unity in $\bbC$, then $\omega I \in \cU(\hilb) \setminus \symthree(\hilb)$.
\end{thm}

\bigskip


When $\dim\, \hilb = d < \infty$, it is easy to see that $J \in \symone(\hilb)$ satisfies $\det\, J \in \{ -1, 1\}$, and thus $V \in \symk(\hilb)$ also satisfies $\det\, V \in \{-1, 1\}$.  The second author of the present paper showed that this is the only obstruction for a unitary operator to lie in $\symfour(\hilb)$.

\begin{thm} \cite[Theorem 3]{Radjavi1969.01} \label{Radjavi1969}
Let $d \in \bbN$ and $V \in \cU(\bbC^d)$.   The following statements are equivalent:
\begin{enumerate}
	\item[(a)]
	$\det\, V \in \{ -1, 1\}$.
	\item[(b)]
	$V = J_1 J_2 J_3 J_4$, where $J_1 J_2$ and $J_3 J_4 \in \symtwo(\bbC^d)$ commute with one another.
\end{enumerate}
\end{thm}
	

\subsection{} \label{sec1.05}
Over the years, a number of people have extended these results to more general settings.   For example, Fillmore~\cite{Fillmore1966} showed that if $\cM$ is a properly infinite von Neumann algebra, then every unitary operator in $\cM$ is a product of at most four symmetries in $\cM$.
Thakare and Baliga~\cite{ThakareBaliga} extended Fillmore's theorem to  properly infinite $AW^*$-algebras.
In the case of a type $II_1$ factor von Neumann algebra $\cM$, it was shown by Broise~\cite{Broise1967} that every unitary element of $\cM$ is a product of finitely many symmetries, while Dowerk and Thom demonstrated that sixteen factors are sufficient for any such unitary.  Recent work of Bhat, Nayak and Shankar~\cite[Theorem~3.9]{BhatNayakShankar2022} demonstrates that every unitary in a type $II_1$ von Neumann algebra $\cM$ (factor or not) is a product of at most six symmetries from $\cM$.     In the same paper, they also prove that the set of products of three symmetries from an arbitrary von Neumann algebra $\cM$ is not norm-dense in the set of unitaries in $\cM$.

In the setting of $C^*$-algebras, Zhang~\cite{ZhangS1993} has shown (amongst many other things) that if $A$ is a unital $C^*$-algebra with real rank zero and topological stable rank one, and if $n \ge 2$ is an integer, then every unitary element of $\bbM_n(A)$ can be factored as a product of at most $\dss 2 \lceil \frac{\ln n}{\ln 2} \rceil$ symmetries from $\bbM_n(A)$ and a diagonal unitary. For a unital purely infinite simple $C^*$-algebra $A$, Leen \cite{Leen1997} proved that every unitary in the identity component of the unitary group of $A$ is a product of finitely many symmetries in $A$. In the case where $A=\bbM_n(C(X))$ (and where $X$ is a finite compact metric space), Phillips~\cite[Theorem 3.9]{Phillips1996} characterized which unitaries in $A$ could be written as a finite product of symmetries in $A$.


\subsection{} \label{sec1.06}
Of course,  a unitary operator $V$ lies in $\symone(\hilb)$ if and only if $\sigma(V) \in \{ -1, 1\}$, or equivalently, if $\sigma(V) \subseteq \bbR$.  Davis~\cite[Theorem 6.3]{Davis1958}, Radjavi and Williams~\cite[Theorem 3]{RadjaviWilliams1969} characterised $\symtwo(\mathcal{H})$ as those unitary operators $V$ which are similar to their adjoints. If $\mathcal{H}=\mathbb{C}^d$, Radjavi and Williams~\cite[Theorem 1]{RadjaviWilliams1969} showed that this is equivalent to the condition that the matrix for $V$ is unitarily equivalent to a matrix with real entries, and Carlson~\cite[Theorem 2]{Carlson1965} showed that this happens if and only if $V$ is a product of two hermitian matrices.


\subsection{} \label{sec1.07}
Whereas the characterisations of $\symone(\hilb),\symtwo(\hilb)$ and $\symfour(\hilb)$ are entirely satisfactory, no characterisation of $\symthree(\hilb)$ currently exists.   The results of the current paper suggest that such a characterisation will likely be very complicated  if it can indeed be determined.  In Section~Three below, we analyse a class of unitary operators $V \in \symthree(\bbC^d)$ whose spectrum consists of exactly two non-real eigenvalues.



\subsection{} \label{sec1.08}
We are also interested in ``simplifying" the unitary orbits of operators.   Recall that the unitary orbit of a given operator $T$ is the set of its unitary conjugates:
\[
\cU(T) = \{ U^* T U: U \in \cU(\hilb) \}. \]
We are interested in the question of whether or not it is possible to use a smaller set, e.g. $\symthree(\hilb)$ or $\symtwo(\hilb)$, instead of $\cU(\hilb)$ in order to generate $\cU(T)$.


%
%
%
%
%
%


\bigskip

\section{A standard form for symmetries relative to a decomposition of the space} \label{section2}


\subsection{} \label{sec2.01}
In the next section, we shall investigate the question of which unitary operators with two distinct eigenvalues may be written as the product of three symmetries.   In order to do this, we shall need to know more about the matrix form of a symmetry relative to certain decompositions of the underlying Hilbert space.  We recall that the \textbf{point spectrum}  $\sigma_p(T)$ of an operator $T \in \bofh$ is the set of eigenvalues of $T$.


\begin{prop} \label{prop2.02}
Let $\hilb$ be a Hilbert space and suppose $\hilb = \hilb_1 \oplus \hilb_2$.  Suppose that $J \in \bofh$ is a symmetry, and that relative to this decomposition,
\[
J = \begin{bmatrix} J_1 & J_2 \\ J_2^* & J_4 \end{bmatrix}. \]
If $1 \not \in \sigma_p(J_1^2) \cup \sigma_p(J_4^2)$, then $\dim\, \hilb_1 = \dim\, \hilb_2$, and after identifying these two spaces, there exists a unitary operator $V = V_0 \oplus I \in \cB(\hilb_1 \oplus \hilb_2)$ such that
\[
V^* J V = \begin{bmatrix} M & N \\ N & - M \end{bmatrix}, \]
where $M= M^*$, and $N = (I-M^2)^{\frac{1}{2}}$ is injective with dense range.
\end{prop}

\begin{pf}
Since $J^2 = I$, we see that $J_1^2 + J_2 J_2^* = I_{\hilb_1}$, or equivalently, $J_2 J_2^* = I_{\hilb_1} - J_1^2$.    Now $\norm J_1 \norm \le 1$ and $1 \not \in \sigma_p(J_1^2)$, so that $\ker\, J_2 J_2^* = \ker (I_{\hilb_1} - J_1^2) = \{ 0\}$.  Hence $J_2^*$ is injective, or equivalently $J_2$ has dense range.

Similarly, $J_2^* J_2 =  I_{\hilb_2} - J_4^2$, and $1 \not \in \sigma_p(J_4^2)$, so that $\ker\, J_2 = \{ 0\}$.    Thus $J_2$ is injective with dense range.    From this it follows that $\dim\, \hilb_1 = \dim\, \hilb_2$.   Let us identify these two spaces and think of $J = \begin{bmatrix} J_1 & J_2 \\ J_2^* & J_4 \end{bmatrix} \in \bbM_2(\cB(\cK))$, where $\cK \simeq\hilb_1 \simeq \hilb_2$.

Let $J_2 = V_0 |J_2|$ now denote the polar decomposition of $J_2$, and observe that the fact that $J_2$ has trivial kernel and dense range implies that $V_0$ is unitary.   Define $M := V_0^* J_1 V_0 = M^*$, $N := |J_2| \ge 0$.  Set $V := V_0 \oplus I$, so that $V \in \bbM_2(\cB(\cK))$ is also unitary.   Then
\[
K:= V^* J V = \begin{bmatrix} M & N \\ N & J_4 \end{bmatrix}. \]
But $K$ is clearly a symmetry (since $J$ is a symmetry and $V$ is unitary), which leads to the equations
\begin{enumerate}
	\item[(eq1)]
		$M^2 + N^2 = I$; and
	\item[(eq2)]
		$M N + N J_4 = 0$.
\end{enumerate}		
Since $N \ge 0$ and since by (eq1), $N^2 = (I-M^2)$, we see that $N$ is the unique positive square root of $(I-M^2)$, i.e. $N= (I-M^2)^{\frac{1}{2}}$.   (Note that $M= M^*$ and that $\norm M \norm \le 1$, and so this is well-defined.). In particular, it follows that $N M = M N$.

From (eq2) we now conclude that $0 = M N + N J_4 = N (M + J_4)$.   Recalling that $\ker\, N = \ker \, J_2 = \{ 0\}$, we conclude that $J_4 = -M$, as required.
\end{pf}


\begin{rem} \label{rem2.03}
With the hypotheses of the above Proposition, suppose now that $d := \dim\, \hilb < \infty$.
Since $M$ is self-adjoint, we may choose a basis for $\cK \simeq \hilb_1 \simeq \hilb_2$ which diagonalises $M$, and which therefore simultaneously diagonalises $N$.    In this case, we find that
\[
J \simeq \oplus_{j=1}^{d/2} \begin{bmatrix} m_j & n_j \\ n_j & -m_j \end{bmatrix}, \]
where $0 < n_j$, $m_j \in (-1, 1)$ and $m_j^2 + n_j^2 = 1$.

It is worth observing that in this setting, $\tr(J) = 0$.
\end{rem}


\begin{prop} \label{prop2.04}
Let $\hilb$ be a Hilbert space and suppose $\hilb = \hilb_1 \oplus \hilb_2$.  Suppose that $J \in \bofh$ is a symmetry, and that relative to this decomposition,
\[
J = \begin{bmatrix} J_1 & J_2 \\ J_2^* & J_4 \end{bmatrix}. \]
Then there exist complementary subspaces $\cM_k, \cN_k$ of $\hilb_k$, $k = 1, 2$,  such that  relative to the decomposition $\hilb = (\cM_1 \oplus \cN_1) \oplus (\cN_2 \oplus \cM_2)$, we have
\[
J = \begin{bmatrix} D_1 & 0 & 0 & 0 \\ 0 & J_{11} & J_{21} & 0 \\ 0 & J_{21}^* & J_{22}  & 0  \\ 0 & 0 & 0 & D_2 \end{bmatrix}, \]
where $D_j = D_j^*$ has spectrum contained in $\{ -1, 1\}$, and $\eta \in \sigma_p(J_{k k})$ implies that $| \eta | < 1$, $k  = 1, 2$.
Moreover, $\dim\, \cN_1 = \dim\, \cN_2$, and by identifying these two spaces with a fixed Hilbert space $\cK$, we see that there exists a unitary operator $V \in \cB(\cK)$ such that relative to the decomposition $\hilb = \cM_1 \oplus \cK \oplus \cK \oplus \cM_2$, we have that
\[
J \simeq D_1 \oplus \begin{bmatrix} M & N \\ N & - M \end{bmatrix} \oplus D_2, \]
where $M = M^*$, and $N = (I-M^2)^{\frac{1}{2}}$ is injective with dense range.
\end{prop}

\begin{pf}
Note that $J_1 = J_1^*$ and $J_4 = J_4^*$.    Let
\begin{itemize}
	\item{}
	$\cM_1 := E(J_1; \{ -1, 1\}) = \mathrm{span} \{\ker\, (J_1 - I), \ker\, (J_1 + I)\}$ and \\
	\item{}
	$\cM_2  := E(J_4; \{ -1, 1\}) = \mathrm{span} \{ \ker\, (J_4 - I), \ker\, (J_4 + I) \}$.
\end{itemize}

Let $\cN_j = \hilb_j \ominus \cM_j$, $j = 1, 2$.   Note that $J$ is a symmetry, and thus $J$ has norm one.   From this it easily follows that $\cM_1$ and $\cM_2$ are reducing for $J$, as well as for $J_1$ and $J_4$ respectively.   In other words, relative to the decomposition $\hilb = \cM_1 \oplus \cN_1 \oplus \cN_2 \oplus \cM_2$, we have
\[
J = \begin{bmatrix} D_1 & 0 & 0 & 0 \\ 0 & J_{1 1} & J_{2 1} & 0 \\ 0 & J_{2 1}^* & J_{2 2}  & 0  \\ 0 & 0 & 0 & D_2 \end{bmatrix}, \]
where  $D_k$ is a symmetry acting on $\cM_k$, $k = 1, 2$, and as such,  $\sigma(D_k) \subseteq \{ -1, 1\}$, $k = 1, 2$.

Since $J$ is a symmetry, so is $J_0 := \begin{bmatrix} J_{1 1} & J_{21} \\ J_{21}^* & J_{2 2} \end{bmatrix}$. Clearly $\sigma_p(J_{k k}) \cap \{ -1, 1\} = \varnothing$ by definition of $\cM_k$, $k = 1, 2$.    Since $J_{k k} = J_{k k}^*$ and $\norm J_{k k} \norm \le 1$,  this implies that $\eta \in \sigma_p(J_{k k})$ implies $\eta \in (-1, 1)$, $k = 1, 2$.

The last statement is now an immediate consequence of Proposition~\ref{prop2.02}.
\end{pf}

We shall refer to the above decomposition of $J$ as the \textbf{standard decomposition} of $J$ relative to $\hilb = \hilb_1 \oplus \hilb_2$.


\begin{rems} \label{rems2.05}
\begin{enumerate}
	\item[(a)]
	Let us now adopt the hypotheses and notation of Proposition~\ref{prop2.04} and further suppose that $d := \dim\, \hilb < \infty$.  Set $r_k := \dim\, \cM_k$, $s_k := \dim\, \cN_k$, $k = 1, 2$.    We may then conclude that there exist orthonormal bases for $\cM_k, \cN_k$, $k = 1, 2$ relative to which
\[
J \simeq \begin{bmatrix} D_1 & 0 & 0 & 0 \\ 0 & M & N & 0 \\ 0 & N & -M  & 0  \\ 0 & 0 & 0 & D_2 \end{bmatrix}, \]
$D_j$ are diagonal symmetries, $M = M^*$ and $0 < N$ are diagonal with $N$ invertible and $M^2 + N^2 = I_{s_1}$.
	\item[(b)]
	In Proposition~\ref{prop2.02} and Remark~\ref{rem2.03}, we observed that if $\hilb = \bbC^s \oplus \bbC^s$ (for some fixed $s \ge 1$) and $J = \begin{bmatrix} J_1 & J_2 \\ J_2^* & J_4 \end{bmatrix} \in \symone(\bbC^{2s})$ is such that $1 \not \in \sigma_p(J_1^2) \cup \sigma_p(J_4^2)$, then there exists $V \in \cU(\bbC^s)$ such that
\[
(V \oplus I_s)^* J (V \oplus I_s) = \begin{bmatrix} M & N \\ N & -M \end{bmatrix}, \]
where $M$ and $N$ are commuting diagonal matrices with $N$ positive and $M^2 + N^2 =I_{s} $.    As such, we concluded that $\tr(J) = 0$.

We now obtain a converse to this result that doesn't require the condition that $1 \not \in \sigma_p(J_1^2) \cup \sigma_p(J_4^2)$.

\end{enumerate}	
\end{rems}


\begin{lem} \label{lem2.06}
Let $s \in \bbN$, $d = 2 s$ and $J \in \symone(\bbC^{d})$.  Suppose that $\tr(J) = 0$.   Then there exist unitary operators $W_1, W_2 \in \cU(\bbC^s)$ such that with $W := W_1 \oplus W_2$,
\[
W^* J W = \begin{bmatrix} M & N \\ N & - M \end{bmatrix}, \]
where $M, N$ are diagonal operators acting on $\bbC^s$.
\end{lem}

\begin{pf}
This follows easily from Proposition~\ref{prop2.02} and Proposition~\ref{prop2.04}.
\end{pf}


\begin{defn} \label{defn2.08}
Let $d \in \bbN$ and $J \in \symone(\bbC^d)$.   We shall say that $J$ is \textbf{balanced} if $\tr (J) \in \{ 0, 1\}$.
\end{defn}


\begin{prop} \label{prop2.07}
Let $d \in \bbN$ and $V \in \symtwo(\bbC^{d})$.  Then there exist $J_1, J_2 \in \symone(\bbC^{d})$  such that $V = J_1 J_2$ and $J_2$ is balanced.

If $W \in  \symthree(\bbC^{d})$, then there exist $J_1, J_2, J_3 \in \symone(\bbC^{d})$ such that $W = J_1 J_2 J_3$, and both $J_2$ and $J_3$ are balanced.
\end{prop}

\begin{pf}
Recall that for any $t \ge 1$, $\symtwo(\bbC^t) = \{ W \in \cU(\bbC^t): W \simeq W^*\}$.   In our case it follows that there exist $0 \le p, q \in \bbZ$ such that $d = 2p + q$ and
\[
V \simeq   \oplus_{j=1}^p \begin{bmatrix} \alpha_j & 0 \\ 0 & \beta_j \end{bmatrix} \oplus \xi I_q, \]
where $\xi \in \{ -1, 1 \}$, and for $1 \le j \le p$, either
\begin{enumerate}
	\item[(i)]
	$\alpha_j \not \in \{ -1, 1\}$ and $\beta_j = \ol{\alpha_j}$, or
	\item[(ii)]
	$\alpha_j \in \{-1, 1\}$ and $\beta_j = - \alpha_j$.
\end{enumerate}	

If $\alpha_j  \not \in \{-1, 1\}$, then  $\beta_j = \ol{\alpha_j}$ and so $\alpha_j \oplus \beta_j = K_j L_j$, where $K_j = \begin{bmatrix} 0 & \alpha_j  \\ \ol{\alpha_j} & 0 \end{bmatrix}$, and $L_j = \begin{bmatrix} 0 & 1 \\ 1 & 0 \end{bmatrix}$.   Clearly $K_j$ and $L_j$ are symmetries, and $\tr(L_j) = 0$.

If $\alpha_j \in \{ -1, 1\}$, then $\beta_j = - \alpha_j$ and so $(\alpha_j \oplus \beta_j) = K_j L_j$, where $K_j = I_2$ and $L_j =  \alpha_i \oplus  \beta_i$.    Again, $K_j$ and $L_j$ are symmetries, and $\tr(L_j) = 0$.

As for $\xi I_q$, we write $q = 2 a +b$, where $b \in \{ 0, 1\}$.    Let $K_{p+1} = \xi  ((I_a \oplus -I_a)  \oplus I_b)$, and let $L_{p+1} = (I_a \oplus -I_a) \oplus I_b$.    Then $\xi I_q = K_{p+1} L_{p+1}$, each of $K_{p+1}$ and $L_{p+1}$ are symmetries, and $\tr(L_{p+1}) = 1$.

We complete this part of the proof by setting $J_1 = \oplus_{j=1}^{p+1} K_j$ and $J_2 = \oplus_{j=1}^{p+1} L_j$.

\bigskip

The statement involving $W \in \symthree(\bbC^d)$ now follows immediately from the above, since by definition of $\symthree(\bbC^d)$, there exist $Q_1 \in \symone(\bbC^d)$ and $V \in \symtwo(\bbC^d)$ such that  $W = Q_1 V$.  By the first half of the proof, we may write $V = Q_2 J_3$, where $Q_2, J_3 \in \symone(\bbC^d)$ and $J_3$ is balanced.   But then $Q_1 Q_2 \in \symtwo(\bbC^d)$, and so from above, we may rewrite $Q_1 Q_2$ as $J_1 J_2$ with $J_1, J_2 \in \symone(\bbC^d)$ where $J_2$  is balanced.  This completes the proof.
\end{pf}


\smallskip

It should be clear that if $W \in \cU(\bbC^d)$ is a product of $n$ symmetries, then $W$ may be written as a product of one symmetry and $n-1$ balanced symmetries.


\section{Unitary elements with small spectrum in $\symthree(\bbC^d)$} \label{section3}


\subsection{} \label{sec3.01}
An ambitious, and at this point seemingly intractable, goal is to try to completely characterise those operators acting upon a finite-dimensional, complex Hilbert space which can be expressed as the product of three symmetries.   Let $2 \le d \in \bbN$.   Since $\symthree(\bbC^d)$ is invariant under unitary conjugation, and since -- up to unitary equivalence -- every unitary operator is completely determined by its eigenvalues and their algebraic multiplicities (or equivalently their geometric multiplicities), the problem reduces to characterising which sets $\{ (\alpha_j, \mu_j)\}_{j=1}^k \subseteq \bbT \times \bbN$ of complex numbers and their multiplicities may appear as the spectrum of an element $V \in \symthree(\bbC^d)$.

We shall obtain a measure of the complexity of this task by considering the following case of the general problem:  for which  $\alpha \ne \beta \in \bbT \setminus \{ -1, 1\}$ can we find $V \in \symthree(\bbC^d)$ such that $\sigma(V) = \{ \alpha, \beta \}$?  To simplify the problem even further, we introduce a condition which -- for lack of a better term -- we refer to as a ``purity" condition.    Condition (ii) below ensures that no direct summand of $V$ can be expressed as a product of three symmetries, while the first condition  is a consequence of the fact that the determinant of a symmetry, and thus of every element of $\symthree(\bbC^d)$, must be real.


\subsection{Condition (P)} \label{sec3.02}
Let $2 \le d \in \bbN$ be an integer and suppose that $d = r+s$ for some $s \le r \in \bbN$.     We say that $\alpha \ne \beta \in \bbT$ satisfy \textbf{Condition (P)} if
\begin{enumerate}
	\item[(i)]
	$\alpha^r \beta^s \in \{ -1, 1 \}$, and
	\item[(ii)]
	if $i, j \in \bbN$ and $\alpha^i \beta^j \in \{ -1, 1\}$, then $i + j \ge d$.
\end{enumerate}
\textbf{Note:} \ \  observe that (ii) implies that $\alpha, \beta \not \in \{ -1, 1\}$.

\bigskip


The proofs of our main results below will require the dimension $d$ of the underlying space to be at least three.    For that reason, we now dispense with the case where $d=2$.  This we can do by citing the following result from~\cite{Radjavi1969.01}.

\begin{prop} \label{prop3.03}
Let $V \in \cU(\bbC^2)$, and suppose that $\sigma(V)= \{ \alpha, \beta\}$ for some $\alpha \ne \beta \in \bbT \setminus \{ -1, 1\}$.  Then $V \in \symthree(\bbC^2)$ if and only if $\det(V) \in \{ -1, 1\}$.
\end{prop}



\subsection{} \label{sec3.04}
We now turn our attention to the case where the dimension $d$ of the underlying Hilbert space is at least $3$.  In this case, the parity of the dimension plays a major role in the description of those unitary operators admitting two distinct, non-real eigenvalues which lie in $\symthree(\bbC^d)$.

Let $\alpha \ne \beta \in \bbT$, $1 \le s \le r \in \bbN$ and $d := r + s$.   Suppose that $V := \alpha I_r \oplus \beta I_s \in \symthree(\bbC^d)$.   By Proposition~\ref{prop2.07}, we may write $V = J_1 J_2 J$, where $J_1, J_2$ and $J$ are symmetries and $J$ is balanced.   As such, the operator $Z := V J = J_1 J_2 \in \symtwo(\bbC^d)$, from which it follows that $Z \simeq Z^*$.    As we shall soon see, we shall be able to decompose $Z$ as a direct sum $Z = Z_0 \oplus Y$, where, for an appropriate choice of $d_0 \in \bbN$, $m_j \in (-1, 1)$, and $n_j \in (0, 1]$ satisfying $m_j^2 + n_j^2 = 1$, $1 \le j \le d_0$, we have
\[
Y \simeq \oplus_{j=1}^{d_0} Y_j, \]
where each $Y_j = \begin{bmatrix} \alpha m_j & \alpha n_j \\ \beta n_j & - \beta m_j \end{bmatrix}$, $1 \le j \le d_0$.   Understanding the spectrum of $Y$, and thus of each $Y_j$, is crucial to our analysis, and so we begin with a few results that apply in both the even-dimensional and the odd-dimensional settings.

%
%

\smallskip

\begin{lem} \label{lem3.06}
Let $m \in (-1, 1)$, $n \in (0,1]$, and suppose that $m^2 + n^2 = 1$.   Let $\alpha \ne \beta \in \bbT$ and  $Y = \begin{bmatrix} \alpha m & \alpha n \\ \beta n & - \beta m \end{bmatrix}$.
	If $\theta \in \sigma(Y)$, then $-\alpha \beta  \ol{\theta} \in \sigma(Y)$.
\end{lem}

\begin{pf}
	An easy computation shows that $Y$ is unitary, whence $\theta \in \bbT$.
	Write $\sigma(Y) = \{ \theta, \eta\}$.   By considering the determinant of $Y$, we obtain:
	\[
	\theta \cdot \eta = (- \alpha \beta), \]
	or equivalently,
	\[
	\eta = (-\alpha \beta) \ol{\theta}. \]
\end{pf}


\begin{lem}\label{lem3.07}
Suppose that $z_j=e^{i \theta_j}\in \bbT$, where $\theta_j\in [0,2\pi)$, $1\leq j\leq 4$. Suppose furthermore that $z_1\neq z_2$ and $z_3\neq z_4$.
The following are equivalent.
\begin{enumerate}
	\item[(a)]
	$(z_4-z_3)/(z_2-z_1) \in \bbR$;
	\item[(b)]
	$(\theta_4+\theta_3)-(\theta_2+\theta_1)\in 2 \pi \bbZ$; and
	\item[(c)]
	$z_3z_4=z_1z_2$.
\end{enumerate}
Moreover,
\[|\frac{z_4-z_3}{z_2-z_1}|<1\text{ \ \ \ if and only if \ \ \  }
|\frac{\sin\frac{\theta_4-\theta_3}{2}}{\sin\frac{\theta_2-\theta_1}{2}}|<1.\]
\end{lem}

\begin{pf}
The equivalence of (b) and (c) is routine, and the proof is omitted.

As to the equivalence of (a) and (b), consider first that by the ``sum-to-product" formulae, for $a, b \in [0, 2 \pi)$,
\begin{align*}
e^{ia}-e^{ib}
	&=\cos a-\cos b+i(\sin a-\sin b)\\
	&= -2\sin \frac{a+b}{2}\sin\frac{a-b}{2}+i2\cos \frac{a+b}{2}\sin\frac{a-b}{2}\\
	&= 2\sin\frac{a-b}{2}(\cos(\frac{a+b}{2}+\frac{\pi}{2})+i\sin(\frac{a+b}{2}+\frac{\pi}{2}) )\\
	&= 2\sin\frac{a-b}{2}e^{i (\frac{a+b}{2}+\frac{\pi}{2})}.
\end{align*}
Therefore,
\begin{align*}
	\frac{z_4-z_3}{z_2-z_1}\in \bbR
		& \text{ if and only if } \frac{e^{i (\frac{\theta_4+\theta_3}{2}+\frac{\pi}{2})}}{e^{i (\frac{\theta_2+\theta_1}{2}+\frac{\pi}{2})}}\in \bbR \\
		& \text{ if and only if } (\frac{\theta_4+\theta_3}{2}+\frac{\pi}{2})-(\frac{\theta_2+\theta_1}{2}+\frac{\pi}{2})\in \pi \bbZ  \\
		& \text{ if and only if }  (\theta_4+\theta_3)-(\theta_2+\theta_1)\in 2 \pi\bbZ.
\end{align*}		
Finally,
\begin{align*}
|\frac{z_4-z_3}{z_2-z_1}|<1
	&  \text{ if and only if }  |\frac{2\sin\frac{\theta_4-\theta_3}{2}}{2\sin\frac{\theta_2-\theta_1}{2}}|<1 \\
	&  \text{ if and only if }  |\frac{\sin\frac{\theta_4-\theta_3}{2}}{\sin\frac{\theta_2-\theta_1}{2}}|<1.
\end{align*}	
\end{pf}


\bigskip

\begin{lem} \label{lem3.08}
Let $3 \le d \in \bbN$ be an integer and suppose that $d = r+s$ for some $s \le r \in \bbN$.     Suppose that $\alpha \ne \beta \in \bbT$ satisfy condition \textbf{\emph{(P)}}.

Let $Z = Z_0 \oplus Y \in \symtwo(\bbC^d)$, and suppose that
\begin{enumerate}
	\item[(\textsc{i})]
	$\{ \alpha \} \subseteq \sigma(Z_0) \subseteq \{ \alpha, -\alpha, \beta, -\beta\}$; and
	\item[(\textsc{ii})]
	$\theta \in \sigma(Y)$ implies that $(-\alpha \beta) \ol{\theta} \in \sigma(Y)$.
\end{enumerate}
Let  $d_0 := \lfloor \frac{d-1}{2} \rfloor$, and for $1 \le j \le d_0$, define the set $\Omega_j(\alpha) := \{ \ol{ (-\alpha \beta)^{j-1} \alpha}, (-\alpha \beta)^j \alpha \}$.  Then

\[
\Omega_1(\alpha) \cup  \Omega_2(\alpha) \cup \cdots \cup \Omega_{d_0}(\alpha)  \subseteq \sigma(Y), \]
where the elements of both sides of this equation are counted according to multiplicity.
\end{lem}	

\begin{pf}
Observe that the condition that $Z \in \symtwo(\bbC^d)$ implies that $Z \simeq Z^*$, and thus
\begin{enumerate}
	\item[(\textsc{iii})]
	$\gamma \in \sigma(Z)$ implies that $\ol{\gamma} \in \sigma(Z)$.
\end{enumerate}


\bigskip
\noindent{\textbf{Step 1.}}

We begin with the observation that if $1 \le j \le d_0$, then
\[
\ol{ (-\alpha \beta)^{j-1} \alpha} \not \in \{ -\alpha, \alpha, -\beta, \beta\}. \]
Indeed, if $\ol{ (-\alpha \beta)^{j-1} \alpha} \in \{-\alpha, \alpha \}$, then
\[
\alpha^{j+1} \beta^{j-1} \in \{ -1, 1\}, \]
whereas if $\ol{ (-\alpha \beta)^{j-1} \alpha} \in \{-\beta, \beta\}$, then
\[
\alpha^j \beta^j \in \{ -1, 1\}. \]
But $(j+1) + (j-1) = j + j = 2j \le 2 d_0 < d$, contradicting condition (ii) of condition \textbf{(P)}.


\bigskip
\noindent{\textbf{Step 2.}}

We now proceed by induction  to show that $\Omega_j(\alpha) \subseteq \sigma(Y)$,  $1 \le j \le d_0$.

\smallskip

First note that $\alpha \in \sigma(Z_0) \subseteq \sigma(Z)$.   By condition $(\textsc{iii})$ listed immediately before \textbf{Step 1} of the proof, $\ol{\alpha} \in \sigma(Z)$.   By \textbf{Step 1}, $\ol{\alpha} = \ol{(-\alpha \beta)^0 \alpha} \not \in \{ -\alpha, \alpha, -\beta, \beta\}$.    Thus $\ol{\alpha} \in \sigma(Y)$.  By condition $(\textsc{ii})$ of the hypotheses, $(-\alpha \beta) \alpha \in \sigma(Y)$, proving that
\[
\Omega_1(\alpha) = \{ \ol{\alpha}, (-\alpha \beta) \alpha\} \subseteq \sigma(Y). \]

Now suppose that $1 \le \kappa \le d_0 -1$ and that $\Omega_\kappa(\alpha) =  \{ \ol{ (-\alpha \beta)^{\kappa-1} \alpha}, (-\alpha \beta)^\kappa \alpha\}   \subseteq \sigma(Y) \subseteq \sigma(Z)$.   By condition $(\textsc{iii})$ above, $\ol{(-\alpha \beta)^\kappa \alpha} \in \sigma(Z)$.  By \textbf{Step 1}, $\ol{(-\alpha \beta)^\kappa \alpha} \not \in \{ -\alpha, \alpha, -\beta, \beta\}$.  Thus $\ol{(-\alpha \beta)^\kappa \alpha} \in \sigma(Y)$.   By condition $(\textsc{ii})$ of the hypotheses,
\[
\Omega_{\kappa + 1}(\alpha) = \{ \ol{(-\alpha \beta)^\kappa \alpha}, (-\alpha \beta)^{\kappa+1} \alpha \} \subseteq \sigma(Y). \]

\smallskip

It therefore follows that $\Omega_1(\alpha) \cup  \Omega_2(\alpha) \cup \cdots \cup \Omega_{d_0}(\alpha) \subseteq \sigma(Y)$, which proves the result.
\end{pf}


\begin{prop} \label{prop3.09}
Let $3 \le d \in \bbN$ be an integer and suppose that $d = r+s$ for some $s \le r \in \bbN$.     Suppose that $\alpha \ne \beta \in \bbT$ satisfy condition \textbf{\emph{(P)}}.

Let $V := \alpha I_r \oplus \beta I_s \in \symthree(\bbC^d)$, $\hilb_1 = \ker\, (V-\alpha I_d)$ and $\hilb_2 = \ker\, (V - \beta I_d)$.   Write $V = J_1 J_2 J$ where $J \in \symone(\bbC^d)$ is balanced, and let
\[
J= \begin{bmatrix} D_1 & 0 & 0 & 0 \\ 0 & M & N & 0 \\ 0 & N & -M & 0 \\ 0 & 0 & 0 & D_2 \end{bmatrix} \]
be the standard decomposition of $J$ relative to $\bbC^d = \hilb_1 \oplus \hilb_2 = (\cM_1 \oplus \cN_1) \oplus (\cN_2 \oplus \cM_2)$.
Finally, suppose that $\dim\, \cM_1 > 0$.
Then, either
\begin{enumerate}
	\item[(a)]
	$d$ is odd, $\dim\, \cM_1 = 1$ and $\dim\, \cM_2 = 0$; or
	\item[(b)]
	$d$ is even and $\dim\, (\cM_1 \oplus \cM_2) =  2$.
\end{enumerate}		
\end{prop}

\noindent{\textbf{Remark.}}  The form of the operator $Z$ appearing in \textbf{Step 1} of this proof will be referred to multiple times in the subsequent results of this section.

\bigskip

\begin{pf}
 As in Lemma~\ref{lem3.08}, we set $d_0 = \lfloor \frac{d-1}{2} \rfloor$.
\bigskip

\noindent{\textbf{Step 1.}}

By Proposition~\ref{prop2.07}, we may write $V = J_1 J_2 J$, where $\textsc{tr}(J) \in \{ 0, 1\}$.   Thus,  $Z := V J = J_1 J_2$ lies in  $\symtwo(\bbC^d)$, and as such, the non-real eigenvalues of $Z$ will appear in conjugate pairs.  That is, the spectrum of $Z$ and $Z^*$ agree, counted according to algebraic multiplicity.

\

Set $r_1 := \dim\, \cM_1 > 0$, $r_2 := \dim\, \cN_1 = \dim\, \cN_2$, and $s_2 := \dim\, \cM_2$.   Relative to this decomposition, we have
\[
V = \begin{bmatrix} \alpha I_{r_1} & 0 & 0 & 0 \\ 0 & \alpha I_{r_2} & 0 & 0 \\ 0 & 0 & \beta I_{r_2} & 0 \\ 0 & 0 & 0 & \beta I_{s_2} \end{bmatrix}, \ \ \ \ \ J = \begin{bmatrix} D_1 & 0 & 0 & 0 \\ 0 & M & N & 0 \\ 0 & N & -M & 0 \\ 0 & 0 & 0 & D_2 \end{bmatrix}. \]
Moreover, $M = \ttt{diag}(m_1, m_2, \ldots, m_{r_2})$ and $N = \ttt{diag}(n_1, n_2, \ldots, n_{r_2})$ are diagonal matrices with $M^2 + N^2 = I_{r_2}$, $-1 < m_j < 1$ for $1 \le j \le r_2$, and $\sigma(D_j) \subseteq \{ -1, 1\}$, $j= 1, 2$.  By multiplying $J$ and $J_2$ by $-1$ if necessary, we may assume without loss of generality that $1 \in \sigma(D_1)$.   It is also not hard to see that by conjugating both $V$ and $J$ by a diagonal unitary of the form $C:= I_{r_1} \oplus C_0 \oplus I_{r_2} \oplus I_{s_1}$, we may assume without loss of generality that $0 < N$.  That is, $C^* V C = V$, while $C^* J C$ has the prescribed form.  \emph{We therefore assume that $V$ and $J$ a priori satisfy these conditions.}

For $1 \le k \le r_2$, set
\[
Y_k = \begin{bmatrix} \alpha m_k & \alpha n_k \\ \beta n_k & -\beta m_k \end{bmatrix}, \]
and define $Z_0 = (\alpha D_1 \oplus \beta D_2)$.  Since $1 \in \sigma(D_1)$, we have that $\alpha \in \sigma(Z_0) = \alpha\, \sigma(D_1) \cup \beta\, \sigma(D_2) \subseteq \{ -\alpha, \alpha, -\beta, \beta\}$.

Thus $Z \simeq Z_0 \oplus Y$, where  $Y = \begin{bmatrix} \alpha M & \alpha N \\ \beta N & - \beta M \end{bmatrix} \simeq \oplus_{k=1}^{r_2} Y_k$.
Moreover, $\sigma (Z) = \sigma (Z_0) \cup \sigma (Y)$, with the eigenvalues counted according to algebraic multiplicity.  If $\theta \in \sigma(Y)$, then $\theta \in \sigma(Y_j)$ for some $1 \le j \le r_2$, and so $(-\alpha \beta) \ol{\theta} \in \sigma(Y_j) \subseteq \sigma(Y)$ by Lemma~\ref{lem3.06}.

By Lemma~\ref{lem3.08}, if $\Omega_j(\alpha) :=  \{ \ol{ (-\alpha \beta)^{j-1} \alpha}, (-\alpha \beta)^j \alpha\}$, $1 \le j \le r_2$, then
\[
\Omega_1(\alpha) \cup  \Omega_2(\alpha) \cup \cdots \cup \Omega_{r_2}(\alpha) \subseteq \sigma(Y). \]

\bigskip

\noindent{\textbf{Step 2.}}

	The idea behind the proof of the Proposition is to prove that this forces $r_2 = d_0$, and to do this, we must show that each $\Omega_j(\alpha)$ occurs separately in the count \emph{according to multiplicity} of the eigenvalues of $Y$.   To do this, we shall apply two separate pairings to the eigenvalues of $Z$.   The first pairing is a consequence of the fact that the non-real eigenvalues of $Z$ must occur in conjugate pairs, counted according to multiplicity.   This is a result of the fact that $Z \simeq Z^*$.   The second pairing arises because of the fact that $Y$ is a direct sum of $2 \times 2$ matrices, each of which has eigenvalues of the form $\{ \theta, (-\alpha \beta) \ol{\theta}\}$ for some $\theta \in \bbC$.

It is worth noting that if $1 \le j \le d_0 - 1$, then $(- \alpha \beta)^j \alpha \not \in \{ -1, 1\}$ by condition \textbf{(P)}, as $2 j + 1 \le 2 (d_0 -1) + 1 < d$.
	
\bigskip

\noindent{\textbf{Step 3.}}

We know that $\alpha \in \sigma(Z)$, and thus $\ol{\alpha} \in \sigma(Z)$.   As we have seen in \textbf{Step 1} of Lemma~\ref{lem3.08}, this forces $\ol{\alpha} \in \sigma(Y) = \cup_{j=1}^{r_2} \sigma(Y_j)$.   By relabelling the indices of the $Y_j$'s if necessary, we may assume without loss of generality that $\ol{\alpha} \in \sigma(Y_1)$.   By Lemma~\ref{lem3.06}, we have that $(-\alpha \beta) \alpha \in \sigma(Y_1)$, and is distinct from $\{ \alpha, \ol{\alpha}\}$ \emph{counted according to multiplicity}.   That is, $\sigma(Y_1) = \Omega_1(\alpha)$, and \emph{counted according to multiplicity}, $\{ \alpha\} \cup \Omega_1(\alpha)$ consists of  three elements of $\sigma(Z)$.   But the non-real eigenvalues of $Z$ occur in conjugate pairs, which implies that the second element of $\Omega_1(\alpha)$, namely $(-\alpha \beta) \alpha$ (unless $d_0 = 1$, in which case $\sigma(Z) = \{ \alpha\} \cup \Omega_1(\alpha)$; as otherwise $1 \le d_0 -1$ - see the last remark of \textbf{Step 2}) must be paired with a \emph{new} eigenvalue $\ol{(-\alpha \beta) \alpha}$ of $Z$ (and \emph{a fortiori} of $Y$).   Without loss of generality, $\ol{(-\alpha \beta) \alpha} \in \sigma(Y_2)$, and by Lemma~\ref{lem3.06}, this forces $\sigma(Y_2) = \Omega_2(\alpha)$.   That is,
\[
\{ \alpha \} \cup \Omega_1(\alpha) \cup \Omega_2(\alpha) \]
consists of five distinct elements of $\sigma(Z)$ when counted according to multiplicity.

The argument now repeats:   the second element of $\Omega_2(\alpha)$, namely $(-\alpha \beta)^2 \alpha \in \sigma(Y_2) \subseteq \sigma(Z)$, and (unless $d_0= 2$), this forces $\ol{(-\alpha \beta)^2 \alpha} \in \sigma(Z)$.   Since the first four elements of $\{ \alpha \} \cup  \Omega_1(\alpha) \cup  \Omega_2(\alpha)$ have already been paired,  $(-\alpha \beta)^2 \alpha \in \sigma(Y)$ must occur as a distinct element of $\sigma(Y)$, which after reindexing if necessary, we may assume comes from $Y_3$.   This in turn forces $\sigma(Y_3) =  \Omega_3(\alpha)$, and so
\[
\{ \alpha \} \cup  \Omega_1(\alpha) \cup  \Omega_2(\alpha) \cup  \Omega_3(\alpha)\]
consists of seven distinct elements of $\sigma(Z)$ when counted according to multiplicity.

By Lemma~\ref{lem3.08}, we see that process continues until we have shown that
\[
\{ \alpha \} \cup\Omega_1(\alpha) \cup  \Omega_2(\alpha) \cup \cdots \cup \Omega_{d_0}(\alpha)  \]
consists of $2 d_0 + 1$ distinct elements of $\sigma(Z)$ when counted according to multiplicity, and without loss of generality,  $\sigma(Y_j) = \Omega_j(\alpha)$, $1 \le j \le d_0$.

\bigskip

\noindent{\textbf{Step 4.}}

It therefore follows that $r_2 \ge d_0$.   But then
\[
d = r_1 + 2 r_2 + s_2 = (2 r_2) + (r_1 + s_2) \ge 2 d_0 + (r_1 + s_2). \]

\begin{itemize}
	\item{}
	If $d$ is odd, then $\dss d_ 0 = \frac{d-1}{2}$, which yields $d \ge 2 (\frac{d-1}{2}) + (r_1 + s_2) = d - 1 + (r_1 + s_2)$.   Given that $r_1 \ge 1$, we conclude that $s_2 = 0$.   Stated otherwise, in this case, $\dim\, \cM_1 = 1$ and $\dim\, \cM_2 = 0$.
	\item{}
	If $d$ is even, then $\dss d_0 = \frac{d-2}{2}$, which yields $d \ge (d-2) + (r_1 + s_2)$.   Given that $r_1 \ge 1$, we either have $\dim\, \cM_1 = r_1 = 2$ and $\dim\, \cM_2 = s_2 = 0$; or $\dim\, \cM_1 = r_1 = 1 = s_2 = \dim\, \cM_2$.   Either way, since $\cM_1$ and $\cM_2$ are orthogonal, $\dim\, (\cM_1 \oplus \cM_2) = 2$.
\end{itemize}
\end{pf}


\begin{rem} \label{rem3.10}
We may summarise the situation as follows:

\bigskip

Let $3 \le d \in \bbN$ be an integer and suppose that $d = r+s$ for some $s \le r \in \bbN$.     Suppose that $\alpha \ne \beta \in \bbT$ satisfy condition \textbf{(P)}.

If $V := \alpha I_r \oplus \beta I_s \in \symthree(\bbC^d)$, then without loss of generality, one of four things takes place, namely:
\begin{enumerate}
	\item[(a)]
	$d$ is odd, $r - s = 1$, $D_1 = I_1$ acts on a space of dimension $1$, and $D_2$ is absent.  That is, $\dim\, \cM_1 = 1$, and $\dim\, \cM_2 = 0$.
	\item[(b)]
	$d$ is even,  $r - s = 2$, $D_1 \simeq I_1 \oplus -I_1 \in \symone(\bbC^2)$, and $D_2$ is absent.
	\item[(c)]
	$d$ is even, $r = s$, $D_1= I_1$, and $D_2 = - I_1$.  That is, $\dim\, \cM_1 = 1 = \dim\, \cM_2$.
	\item[(d)]
	$d$ is even, and $r = s$, and $\dim\, \cM_1 = \dim\, \cM_2 = 0$.
\end{enumerate}		
In cases (b) and (c), the fact that $\textsc{tr}(D_1 \oplus D_2) = 0$ is a consequence of the fact that $\textsc{tr}(J) = 0$.   As we shall see, each of these four scenarios is possible, but each also places further restrictions upon the choice of $\alpha$ and $\beta$.  But first, we require a technical Remark.
\end{rem}


\begin{rem} \label{rem3.11}
Let $d_0 \in \bbN$, and suppose that $\alpha, \beta \in \bbT$ satisfy $\alpha^{d_0 +1} \beta^{d_0} \in \{ -1, 1\}$.   We claim that the condition
\[
| 1 - \alpha^{2j +1} \beta^{2j -1}| < | \alpha - \beta |, \ \ \ 1 \le j \le d_0 \tag{$\ast$} \]
is equivalent to the condition
\[
| 1 - (\alpha \beta)^{2 j - 1}| < | \alpha - \beta |, \ \ \ 1 \le j \le d_0 \tag{$\ast \ast$}. \]

Indeed, fixing $1 \le j \le d_0$, we find that
\begin{align*}
|1 - \alpha^{2 j+1} \beta^{2 j - 1} |
	&= |1 - \alpha^{-2j -1} \beta^{-2j +1} | \\
	&= | 1 - \alpha^{2d_0 - 2j + 1} \beta^{2d_0 - 2j +1} |\\
	&= | 1 - \alpha^{2 k - 1} \beta^{2 k - 1}| \\
\end{align*}
where $k := d_0 - j + 1$.   Since  $\{ d_0 - j + 1: 1 \le j \le d_0\} = \{ 1, 2, \ldots, d_0\}$, we see that $(\ast)$ and $(\ast \ast)$ are equivalent.	
\end{rem}


\begin{center}
\textsc{Case One:  The odd-dimensional setting}
\end{center}

\bigskip


\begin{thm} \label{thm3.12}
Let $3 \le d \in \bbN$ be an odd integer and suppose that $d = r+s$ for some $s < r \in \bbN$.   Suppose that  $\alpha \ne \beta \in \bbT$ satisfy condition~$\textbf{\emph{(P)}}$, and set  $d_0 = \dss \frac{d - 1}{2}$.

Let  $V := \alpha I_r \oplus \beta I_s \in \cU(\bbC^d)$.  The following statements are equivalent:
\begin{enumerate}
	\item[(a)]
	$V \in \symthree(\bbC^d)$.
	\item[(b)]
	$r = d_0 + 1$, $s =d_0$,  and
	\[
	| 1 - (\alpha \beta)^{2 j - 1}|   < | \alpha - \beta|, 1 \le j \le d_0. \]	
\end{enumerate}			
\end{thm}

\begin{pf}

\noindent{(a) implies (b)}.
\ \ \
By Proposition~\ref{prop3.09}, we see that $V = J_1 J_2 J \in \symthree(\bbC^d)$ implies that $r = d_0 +1$ and $s = d_0$.  Furthermore, without loss of generality
\[
Z := V J \simeq \alpha I_1 \oplus (\oplus_{j=1}^{d_0} Y_j), \]
where each $Y_j$ is of the form $\begin{bmatrix} \alpha m_j & \alpha n_j \\ \beta n_j & -\beta m_j \end{bmatrix}$ as in \textbf{Step 1} of the proof of Proposition~\ref{prop3.09}, and  $\sigma(Y_j) = \Omega_j(\alpha) := \{ \ol{(- \alpha \beta)^{j-1} \alpha}, (- \alpha \beta)^j \alpha \}$, $1 \le j \le d_0$.

\bigskip

Next, for each $1 \le j \le d_0$, we have
\[
\tr(Y_j)  = \ol{(- \alpha \beta)^{j-1} \alpha} + (- \alpha \beta)^j \alpha, \]
and thus
\[
|1 - \alpha^{2j+1} \beta^{2j-1} | = |\tr(Y_j)| = |(\alpha  - \beta) m_j|  < | \alpha - \beta|, 1 \le j \le d_0. \]

By Remark~\ref{rem3.11}, this is equivalent to the condition
\[
|1 - (\alpha \beta)^{2j-1} | < | \alpha - \beta|, 1 \le j \le d_0. \]

\bigskip

\noindent{(b) implies (a)}.
\ \ \

The proof of the sufficiency of condition (b) to imply that $V = \alpha I_{d_0 +1}  \oplus \beta I_{d_0} \in \symthree(\bbC^d)$ is mostly a matter of reverse-engineering the  argument of Proposition~\ref{prop3.09}, along with an application of Lemma~\ref{lem3.07}.

\bigskip

To that end, suppose that (b) holds, and for $1 \le j \le d_0$, consider
\[
m_j := \frac{\ol{(-1)^{j-1} \alpha^{j} \beta^{j-1}} + (-1)^j \alpha^{j+1} \beta^j}{\alpha - \beta} = (-1)^{j-1}  \frac{\ol{\alpha^j \beta^{j-1}} - \alpha^{j+1} \beta^j}{\alpha - \beta}. \]
It follows immediately from Lemma~\ref{lem3.07} that each $m_j \in \bbR$, and the hypotheses of (b) imply that $|m_j| < 1$ for all $1 \le j \le d_0$.   Let $n_j := \sqrt{1-m_j^2}$, so that $0 < n_j$ and $m_j^2 + n_j^2 = 1$, $1 \le j \le d_0$.

Set $Y_j := \begin{bmatrix} \alpha m_j & \alpha n_j \\ \beta n_j & - \beta m_j \end{bmatrix}$, $1 \le j \le d_0$ and $Y = \oplus_{j=1}^{d_0} Y_j$.      Then an application of Lemma~\ref{lem3.06} shows that $\sigma(Y_j) =\Omega_j(\alpha)$, $1 \le j \le d_0$.  Since $\alpha^{d_0+1} \beta^{d_0} \in \ \{ -1,1\}$, we see that $\sigma(Y_{d_0}) \cap \{ -1, 1\} \ne \varnothing$, and hence if we  define $Z = \alpha I_1 \oplus Y$, then $\sigma(Z) = \sigma(Z^*)$, \emph{counted according to multiplicity}.

Since $\alpha I_1$ and $Y$ are easily seen to be unitary, it therefore follows that $Z \in \symtwo(\bbC^d)$, and with $J$ defined as in \textbf{Step~1} of the proof of Proposition~\ref{prop3.09}, we see that $J \in \symone(\bbC^d)$, from which we conclude that $Z J = \alpha I_{d_0 + 1} \oplus \beta I_{d_0} \in \symthree(\bbC^d)$.
\end{pf}


\begin{eg} \label{eg3.13}
Theorem~\ref{thm3.12} allows us to produce explicit examples of  unitary operators in $\symthree(\bbC^d)$ which have only two points $\alpha, \beta$ in their spectrum, where $\alpha, \beta \in \bbT$ satisfy condition \textbf{(P)}.

Suppose that $d_0 = 2$, and $d=2 d_0 +1=5$.  Let
\[\alpha=e^{i \frac{58}{360} 2\pi},~~\beta=e^{i \frac{183}{360}2\pi}.\] Then $\alpha^3\beta^2=-1$,
and $\alpha^i \beta^j \in \{ -1, 1\}$ implies that $i+j \ge d$. We claim that
\[
V =\alpha I_3\oplus \beta I_2\in \symthree(\bbC^5). \]
By Theorem~\ref{thm3.12}, it suffices to check that
\[
|\frac{\alpha\beta-1}{\beta-\alpha}|<1 \text{ and } ~~|\frac{(\alpha\beta)^3-1}{\beta-\alpha}|<1.\]
By Lemma~\ref{lem3.07}, this in turn reduces to the problem of verifying that
\[
|\frac{\sin\frac{\arg{\alpha\beta}-0}{2}}{\sin\frac{\arg{\beta}-\arg{\alpha}}{2}}|<1 \text{ and }~~|\frac{\sin\frac{\arg{(\alpha\beta)^3}-0}{2}}{\sin\frac{\arg{\beta}-\arg{\alpha}}{2}}|<1,\]
where $\arg(z)$ is the main value of the argument of the complex number $z$,  which is the unique number in $[0,2\pi)$ such that $e^{i\arg(z)}=z$ for $z\neq 0$.

Note that
\[\frac{\arg{(\alpha\beta)}-0}{2}=\frac{120.5}{360}2\pi,\ \ \ \ \ \frac{\arg{(\alpha\beta)^3}-0}{2}=\frac{1.5}{360}2\pi, \ \ \ \ \ \frac{\arg{\beta}-\arg{\alpha}}{2}=\frac{62.5}{360}2\pi.\]

Thus, by the elementary fact that $\sin \theta$ is an increasing function on $[0, \frac{\pi}{2}]$,
\begin{align*}
|\sin\frac{\arg{\alpha\beta}-0}{2}|
	&=\sin \frac{120.5}{360}2\pi=\sin (\pi-\frac{120.5}{360}2\pi) \\
	&=\sin\frac{59.5}{360}2\pi \\
	&<\sin \frac{62.5}{360}2\pi=|\sin \frac{62.5}{360}2\pi|,
\end{align*}	
and
\[|\sin\frac{\arg{(\alpha\beta)^3}-0}{2}|=\sin\frac{1.5}{360}2\pi<\sin \frac{62.5}{360}2\pi=|\sin \frac{62.5}{360}2\pi|.\]
By Theorem~\ref{thm3.12}, we conclude that $V \in \symthree(\bbC^5)$.
\end{eg}


\vskip 1 cm

\begin{center}
\textsc{Case Two:  The even-dimensional setting,  $\dim\, \cM_1 = 2$ and $\dim\, \cM_2 = 0$.}
\end{center}


\bigskip

The following proof is very similar to that of the previous case.

\begin{thm} \label{thm3.14}
Let $4 \le d \in \bbN$ be an even integer and suppose that $d= r+ s$ for some $s < r \in \bbN$.   Suppose that  $\alpha \ne \beta \in \bbT$ satisfy condition~$\textbf{\emph{(P)}}$.   Let $d_0 = \frac{d-2}{2}$, and $V := \alpha I_r \oplus \beta I_s \in \cU(\bbC^d)$.  The following statements are equivalent:
\begin{enumerate}
	\item[(a)]
	$V \in \symthree(\bbC^d)$.
	\item[(b)]
	$r = d_0 +2$, $s = d_0$, $\alpha^r \beta^s = (-1)^{d_0 + 1}$, and
	\[
	|1 - \alpha^{2j+1} \beta^{2j-1} |  < | \alpha - \beta|, \ \ 1 \le j \le  d_0. \]	
\end{enumerate}			
\end{thm}

\begin{pf}

\noindent{(a) implies (b)}.
\ \ \
By Remark~\ref{rem3.10} (b), we see that $V  \in \symthree(\bbC^d)$ implies that $r = d_0 + 2$ and $s = d_0$.  Furthermore, writing $V = J_1 J_2 J$, where $\tr(J) = 0$, we find that without loss of generality
\[
Z := V J = \alpha D_1 \oplus (\oplus_{j=1}^{d_0} Y_j), \]
where $D_1 \simeq I_1 \oplus -I_1 \in \symone(\bbC^2)$, and $\sigma(Y_j) = \Omega_j(\alpha)  := \{ \ol{ (-\alpha \beta)^{j-1} \alpha}, (-\alpha \beta)^j \alpha\}$, $1 \le j \le d_0$.   As before, each $Y_j$ is of the form $\begin{bmatrix} \alpha m_j & \alpha n_j \\ \beta n_j & -\alpha m_j \end{bmatrix}$, where $m_j \in (-1, 1)$ and $n_j = \sqrt{1-m_j^2}$, $1 \le j \le d_0$.

\bigskip

Next, for each $1 \le j \le d_0$, we have
\[
\tr(Y_j) =  \ol{ (-\alpha \beta)^{j-1} \alpha} + (-\alpha \beta)^j \alpha, \]
and thus
\[
|1 - \alpha^{2j+1} \beta^{2j-1} | = |\tr(Y_j)| = |(\alpha  - \beta) m_j|  < | \alpha - \beta|, 1 \le j \le d_0. \]

The fact that the eigenvalues of $Z$ come in conjugate pairs also implies that
\[
-\alpha = \ol{ (-1)^{d_0} \alpha^{d_0+1} \beta^{d_0 }}, \]
so that $\alpha^{d_0 + 2} \beta^{d_0}  = (-1)^{d_0 + 1}$.

\bigskip
\noindent{(b) implies (a)}.
\ \ \

As in the previous case, the proof of the sufficiency of conditions (b) to imply that $V = \alpha I_{d_0 +1}  \oplus \beta I_{d_0 - 1} \in \symthree(\bbC^d)$ is mostly a matter of reverse-engineering the  argument of Proposition~\ref{prop3.09}, along with an application of Lemma~\ref{lem3.07}.

\bigskip

To that end, suppose that (b) holds, and for $1 \le j \le d_0$, consider
\[
m_j := \frac{\ol{(-1)^{j-1} \alpha^{j} \beta^{j-1}} + (-1)^j \alpha^{j+1} \beta^j}{\alpha - \beta} = (-1)^{j-1}  \frac{\ol{\alpha^j \beta^{j-1}} - \alpha^{j+1} \beta^j}{\alpha - \beta}. \]
It follows immediately from Lemma~\ref{lem3.07} that each $m_j \in \bbR$, and the hypotheses of (b) imply that $|m_j| < 1$ for all $1 \le j \le d_0$.   Let $n_j := \sqrt{1-m_j^2}$, so that $0 < n_j$ and $m_j^2 + n_j^2 = 1$, $1 \le j \le d_0$.

Set $Y_j := \begin{bmatrix} \alpha m_j & \alpha n_j \\ \beta n_j & - \beta m_j \end{bmatrix}$, $1 \le j \le d_0$ and $Y = \oplus_{j=1}^{d_0} Y_j$.      Then an application of Lemma~\ref{lem3.06} shows that $\sigma(Y_j) = \Omega_j(\alpha)$, $1 \le j \le d_0$.   Set $D_1 = I_1 \oplus -I_1$, so that $\sigma(\alpha D_1) = \{ \alpha, -\alpha\}$, and finally set $Z := \alpha D_1 \oplus Y$.

The above construction allows us to pair the first eigenvalue (namely $\alpha$) coming from $\alpha D_1$ with the first eigenvalue of $Y_1$, and to pair  the second eigenvalue of $Y_{j}$ with the first eigenvalue of $Y_{j+1}$, $1 \le j \le d_0 - 1$   to obtain $d_0$ conjugate pairs of eigenvalues of $Z$.

The remaining pair of eigenvalues  of $Z$ are $-\alpha$ (coming from $\alpha D_1$) and $(-1)^{d_0} \alpha^{d_0+1} \beta^{d_0}$, the second eigenvalue of $Y_{d_0}$.  But the hypothesis of (b)  that $\alpha^{d_0 + 2} \beta^{d_0} = (-1)^{d_0 +1}$ implies that
\[
(-1)^{d_0} \alpha^{d_0+1} \beta^{d_0} = - \ol{\alpha}. \]

Thus the eigenvalues of $Z$ appear in conjugate pairs, so that $Z \in \symtwo(\bbC^d)$.   Defining $J$ as in Proposition~\ref{prop3.09}, we find that $V := Z J \in \symthree(\bbC^d)$.
\end{pf}


\begin{eg} \label{eg3.15}
Let us consider the setting where $d = 8$, $r = 5$, $s = 3$, and $d_0 = 3$.

\bigskip

Let $\alpha = e^{2\pi i/16}$, and set $\beta = - \alpha$.   Observe that $\alpha^5 \beta^3 = - \alpha^8 = 1 = (-1)^{d_0 + 1}$.   If $j+k < 8$, then $\alpha^j \beta^k = \pm e^{ (j+k) 2 \pi i/16} \not \in \{ -1, 1\}$.
In particular, $|1 - \alpha^{2j+1} \beta^{2j-1} |  <2=| \alpha - \beta|, \ \ 1 \le j \le  d_0.$
Let $V = \alpha I_5 \oplus \beta I_3 \in \cU(\bbC^8)$, and note that  $\det(V)=1$.  By Theorem \ref{thm3.14}, $V\in \symthree(\bbC^8)$.
\end{eg}


\bigskip

\begin{center}
\textsc{Case Three:  The even-dimensional setting, $\dim\, \cM_1 = 1 = \dim\, \cM_2$.}
\end{center}

\bigskip


\begin{rem} \label{rem3.16}
We now consider case (c) of Remark~\ref{rem3.10}, so that $d$ is even, $r = s $ and $\dim\, \cM_1 = 1 = \dim\, \cM_2$.  We begin by observing that this case only occurs if the dimension $d$ of the underlying space is not divisible by $4$.

\bigskip

Indeed, in case (c) of Remark~\ref{rem3.10} we have that $D_1 = I_1 = -D_2$.    Maintaining the same notation we have been using above,  we find that
\[
Z := VJ = (\alpha I_1 \oplus -\beta I_1) \oplus Y, \]
where $Y = \oplus_{j=1}^{d_0} Y_j$, and  $d_0 := \frac{d}{2} - 1$.   The only way that we can have $Z \simeq Z^*$ is if the non-real eigenvalues of $Z$  appear in conjugate pairs, in which case the second eigenvalue of $Y_{d_0}$, namely $(-\alpha \beta)^{d_0} \alpha$, is paired with $-\beta$, or to be more precise,
\[
(- \alpha \beta)^{d_0} \alpha = \ol{-\beta}. \]
If $d$ is divisible by $4$, then $d_0$ is odd, and so $(-1)^{d_0} = -1$.
This implies that
\[
\alpha^{d/2} \beta^{d/2} = \alpha^{d_0 + 1} \beta^{d_0 + 1} = 1, \]
and thus
\[
(\alpha)^{d/4} (\beta)^{d/4} \in \{-1, 1\}, \]
contradicting condition \textbf{(P)}.
\end{rem}


\bigskip

\begin{thm} \label{thm3.17}
Let $4 \le d \in \bbN$ be an even integer and suppose that $d = r+s$ for some $s \le r \in \bbN$.   Suppose that  $\alpha \ne \beta \in \bbT$ satisfy condition $\textbf{\emph{(P)}}$.  Let  $V := \alpha I_r \oplus \beta I_s.$
The following statements are equivalent:
\begin{enumerate}
	\item[(a)]
	$V \in \symthree(\bbC^d)$, and $\dim\, \cM_1 = 1$ in the standard decomposition of $J$ from Proposition~\ref{prop3.09}.
	\item[(b)]
	$d$ is not divisible by $4$, $s = d/2 =r$, $\alpha^{d/2} \beta^{d/2} = -1$, and with $d_0 = d/2 - 1$,
	\[
	| 1 - \alpha^{2 j+1}  \beta^{2 j - 1}|   < | \alpha - \beta|, \ \ \ \ \ 1 \le j \le d_0. \]	
\end{enumerate}			
\end{thm}

\begin{pf}
\begin{enumerate}
	\item[(a)] implies (b).
	By Proposition~\ref{prop3.09} and Remark~\ref{rem3.10}, we must have $r = s = d/2$.  Moreover,
	\[
	Z = (\alpha I_1 \oplus -\beta I_1) \oplus Y, \]
	where $Y = \oplus_{j=1}^{d_0} Y_j$ and $\sigma(Y_j) = \Omega_j(\alpha)$, $1 \le j \le d_0$.
	
	In order for the non-real eigenvalues of $Z$ to appear in conjugate pairs (which is necessary as $Z \in \symtwo(\bbC^d)$, the second eigenvalue of $Y_{d_0}$, namely $(-\alpha \beta)^{d_0} \alpha$ must be paired with the second eigenvalue of $Z_0$, namely $-\beta$.   More precisely,
	\[
	(-\alpha \beta)^{d_0} \alpha = \ol{-\beta}. \]
	From this and the fact that $d_0$ is
even, it follows that
	\[
	\alpha^{d/2} \beta^{d/2} = \alpha^{d_0+1} \beta^{d_0+1} = -1. \]
	
	Finally, the trace condition on each $Y_j$ always holds, namely:
	\[
	\ \ \ |\textsc{tr}(Y_j)| = | \ol{ (-\alpha \beta)^{j-1} \alpha} + (-\alpha \beta)^{j} \alpha |  = | (\alpha - \beta) m_j| < | \alpha - \beta|, \ \ \ \ \ 1 \le j \le d_0, \]
	which, as we have seen more than once, is equivalent to the inequality listed in (b) above.
	\item[(b)] implies (a).
	
	Set $Z_0 = \alpha I_1 \oplus -\beta I_1$, and for $1 \le j \le d_0$, choose $Y_j = \begin{bmatrix} \alpha m_j & \alpha n_j \\ \beta n_j & - \beta m_j \end{bmatrix}$ such that the eigenvalues of $Y_j$ are $\Omega_j(\alpha) = \{ \ol{(-\alpha \beta)^{j-1} \alpha}, (- \alpha \beta)^j \alpha\}$, $1 \le j \le d_0$.  Let $Z = Z_0 \oplus Y$, where $Y = \oplus_{j=1}^{d_0} Y_j$.   The condition $\alpha^{d/2} \beta^{d/2} = -1$ implies that the second eigenvalue of $Y_{d_0}$ is paired with the second eigenvalue of $Z_0$, whereas the first eigenvalue of $Z_0$ is paired with the first eigenvalue of $Y_1$ and the second eigenvalue of $Y_j$ is paired with the first eigenvalue of $Y_{j+1}$, $1 \le j \le d_0 -1$.
	
	In other words, $Z \simeq Z^*$ and therefore $Z \in \symtwo(\bbC^d)$.   With $J = \begin{bmatrix} 1 & 0 & 0 & 0 \\ 0 & M & N & 0 \\ 0 & N & -M & 0 \\ 0 & 0 & 0 &-1 \end{bmatrix}$, we see that $J \in \symone(\bbC^d)$, and $Z J = V \in \symthree(\bbC^d)$.
\end{enumerate}
\end{pf}


\begin{eg} \label{eg3.18}
Let us consider the setting where $d = 10$.

\bigskip

Let $\alpha = e^{2\pi i/20}$, and set $\beta = - \alpha$.   Observe that $\alpha^5 \beta^5 = - \alpha^{10} = -1$.   If $j+k < 10$, then $\alpha^j \beta^k = \pm e^{ (j+k) 2 \pi i/20} \not \in \{ -1, 1\}$.
Here $d_0 = 4$ as well. In particular, $|1 - \alpha^{2j+1} \beta^{2j-1} |  <2=| \alpha - \beta|, \ \ 1 \le j \le  d_0.$

\bigskip

Let $V = \alpha I_5 \oplus \beta I_5 \in \cU(\bbC^{10})$.   Clearly $\det (V) =\alpha^5\beta^5= -1$.
 Then by Theorem \ref{thm3.17}, $V\in \symthree(\bbC^{10})$.
\end{eg}



\bigskip

\begin{center}
\textsc{Case Four:  The even-dimensional setting,  $\dim\, \cM_1 = 0 = \dim\, \cM_2$.}
\end{center}



\subsection{} \label{sec3.19}
Before embarking on the last case, we shall prove an analogue of Lemma~\ref{lem3.08}.  Given $k \in \bbN$, $\alpha, \beta \in \bbT$ and $\zeta \in \bbC$, we define $\Delta_k(\zeta) := \{ \ol{(-\alpha \beta)^{k-1}} \zeta, (- \alpha \beta)^k \ol{\zeta}\}$.


\begin{lem} \label{lem3.20}
Let $4 \le d \in \bbN$ be an even integer and suppose that $d = r+s$ for some $s \le r \in \bbN$.   Suppose that  $\alpha \ne \beta \in \bbT$ satisfy condition $\textbf{\emph{(P)}}$.

Let $d_0 := \frac{d}{2}$ and suppose that $Y = \oplus_{j=1}^{d_0} Y_j$, where each $Y_j \in \cU(\bbC^2)$ has the property that there exists $\eta_j \in \bbT$  for which $\sigma(Y_j) = \Delta_1(\eta_j)$.   Finally, suppose that $Y \in \symtwo(\bbC^d)$.
\begin{enumerate}
	\item[(a)]
	If $\xi \in \sigma(Y) \cap \{ -1, 1\} \ne \varnothing$, then $\sigma (Y) = \cup_{j=1}^{d_0} \Delta_j(\xi)$, and $(-\alpha \beta)^{d_0} \in \{ -1, 1\}$.
	\item[(b)]
	If $\sigma(Y) \cap \{ -1, 1\} = \varnothing$ and $\theta \in \sigma(Y)$, then $\sigma (Y) = \cup_{j=1}^{d_0} \Delta_j(\theta)$, and $(-\alpha \beta)^{d_0}  = 1$.
\end{enumerate}	
\end{lem}

\begin{pf}
The statement that $Y \in \symtwo(\bbC^d)$ implies that $Y \simeq Y^*$, and thus the non-real eigenvalues of $Y$ must appear in conjugate pairs, counted according to algebraic multiplicity.  In the proof below, the phrase ``without loss of generality" will always be meant to signify ``after a permutation of the indices, if required".

\begin{enumerate}
	\item[(a)]
	Suppose that $\xi \in \sigma(Y) \cap \{ -1, 1\} \ne \varnothing$;  without loss of generality, $\xi \in \sigma(Y_1)$.   Then, by hypothesis, $\sigma(Y_1) = \Delta_1(\xi)$, implying that $(-\alpha \beta) \xi = (-\alpha \beta) \ol{\xi} \in \sigma(Y_1) \subseteq \sigma(Y)$.  Note that $(-\alpha \beta) \xi  \not \in \bbR$, for otherwise $\alpha \beta \in \{ -1, 1\}$, contradicting the condition \textbf{(P)} hypothesis.     Since $(-\alpha \beta) \xi$ must be paired with a non-real eigenvalue of $Y$, we see that  $\ol{ (-\alpha \beta)} \xi \in \sigma(\oplus_{j=2}^{d_0} Y_j)$.   Without loss of generality, $\ol{ (-\alpha \beta)} \xi  \in \sigma(Y_2)$, and thus
	\[
	\sigma(Y_2) = \{ \ol{ (-\alpha \beta)} \xi, (-\alpha \beta)^2 \xi\}  = \Delta_2 (\xi) = \Delta_1(\ol{ (-\alpha \beta)} \xi). \]
	If $\sigma(Y_2) \cap \{ -1, 1\} \ne \varnothing$, then $(-\alpha \beta)^2 \xi \in \{ -1, 1\}$, and so by the condition \textbf{(P)} hypothesis,  $d_0 = 2$ and we are done.  Otherwise, $\sigma(Y_2) \cap \{ -1, 1\} = \varnothing$ and the process continues.
	
	In general, suppose that $1 \le \kappa < d_0$ and that we have shown that $\sigma(\oplus_{j=1}^\kappa Y_j) = \cup_{j=1}^\kappa \Delta_j(\xi)$, and that $\sigma(Y_j) \cap \{ -1, 1\} = \varnothing$, $2 \le j \le \kappa$.    For $1 \le j \le \kappa-1$, the second eigenvalue of $Y_j$, namely $(-\alpha \beta)^{j} \xi$, is paired with the first eigenvalue $\ol{ (-\alpha \beta)}^{j} \xi$ of $Y_{j+1}$.    Since $\sigma(Y_\kappa) \cap \{ -1, 1\} = \varnothing$, the second eigenvalue of $Y_\kappa$, namely $(-\alpha \beta)^\kappa \xi$, cannot stand alone,  and thus must be paired with an eigenvalue of $\sigma( \oplus_{\kappa + 1}^{d_0} Y_j)$.   We conclude that without loss of generality, $\ol{(-\alpha \beta)^{\kappa +1}} \xi \in \sigma(Y_{\kappa +1}) \subseteq \sigma(Y)$, and thus $\sigma(Y_{\kappa +1}) = \Delta_{\kappa+1}(\xi)$.   If $\kappa +1 < d_0$, then
	\[
	\sigma(Y_{\kappa +1}) \cap \{ -1, 1\} = \Delta_{\kappa + 1} (\xi) =   \{ \ol{(-\alpha \beta)^\kappa} \xi, (-\alpha \beta)^{\kappa +1} \xi\} \cap \{ -1, 1\} = \varnothing \]
	follows easily from  condition \textbf{(P)}.
	
	As such, this process repeats $d_0$ times, implying that $\oplus_{j=1}^{d_0} \Delta_j(\xi) = \sigma(Y)$.
	
	Finally, the second eigenvalue $(-\alpha \beta)^j \xi$ of $Y_j$ and thus the first eigenvalue $\ol{(-\alpha \beta)^{j}} \xi$ of $Y_{j+1}$ are not real for $1 \le j \le d_0 - 1$, while the first eigenvalue of $Y_1$, namely $\xi$ is real,  which implies that the second eigenvalue of $Y_{d_0}$ must be real.    In other words,
\[
(-\alpha \beta)^{d_0} \xi \in \{ -1, 1\}. \]
	From this the desired conclusion follows.
	\item[(b)]
	The proof is very similar to that of (a).  Suppose $\sigma(Y) \cap \{ -1, 1\} = \varnothing$ and that $\theta \in \sigma(Y)$.   Without loss of generality, we may assume that $\theta \in \sigma(Y_1)$, and thus $\sigma(Y_1) = \Delta_1(\theta) = \{ \theta, (-\alpha \beta) \ol{\theta}\}$.
	
	Note that $\sigma(Y_1) \cap \{ -1, 1\} = \varnothing$ by hypothesis, and thus $\ol{(-\alpha \beta)} \theta \not \in \sigma(Y_1)$.   Indeed, if $\ol{(-\alpha \beta)} \theta = \theta$, then $\alpha \beta = - 1$, contradicting condition \textbf{(P)}, while $\ol{(-\alpha \beta)} \theta = (-\alpha \beta) \ol{\theta}$ implies that $(-\alpha \beta) \ol{\theta} \in \{ -1, 1\}$,  a contradiction.  Thus $\ol{(-\alpha \beta)} \theta \in \sigma(\oplus_{j=2}^{d_0} Y_j)$.   Without loss of generality, $\ol{(-\alpha \beta)} \theta \in \sigma(Y_2) = \Delta_1(\ol{(-\alpha \beta)} \theta) = \Delta_2(\theta)$.
	
	Since $\sigma(Y) \cap \{ -1, 1\} = \varnothing$ by hypothesis, this means that the second eigenvalue of $Y_2$, namely $(-\alpha \beta)^2 \ol{\theta}$, has yet to be paired and the process continues.
	
		In general, suppose that $1 \le \kappa < d_0$ and that we have shown that $\sigma(\oplus_{j=1}^\kappa Y_j) = \cup_{j=1}^\kappa \Delta_j(\theta)$.     For $1 \le j \le \kappa-1$, the second eigenvalue of $Y_j$, namely $(-\alpha \beta)^{j} \ol{\theta}$, is paired with the first eigenvalue $\ol{ (-\alpha \beta)}^{j} \theta$ of $Y_{j+1}$.    Since $\sigma(Y_\kappa) \cap \{ -1, 1\} = \varnothing$, the second eigenvalue of $Y_\kappa$, namely $(-\alpha \beta)^\kappa \ol{\theta}$, cannot stand alone.   Also, it cannot be paired with the first eigenvalue $\theta$ of $Y_1$, for otherwise $(-\alpha \beta)^\kappa \ol{\theta} = \ol{\theta}$ implies that $(\alpha \beta)^\kappa = 1$, contradicting condition \textbf{(P)}.   Thus $(-\alpha \beta)^\kappa \ol{\theta}$ must be paired with an eigenvalue of $\sigma( \oplus_{\kappa + 1}^{d_0} Y_j)$.   We conclude that without loss of generality, $\ol{(-\alpha \beta)^{\kappa +1}} \theta \in \sigma(Y_{\kappa +1}) \subseteq \sigma(Y)$, and thus $\sigma(Y_{\kappa +1}) = \Delta_{\kappa+1}(\theta)$.
			
	As such, this process repeats $d_0$ times, implying that $\oplus_{j=1}^{d_0} \Delta_j(\xi) = \sigma(Y)$.

	Since none of the eigenvalues of $Y$ are real, and since the only eigenvalues we have yet to pair are the first eigenvalue of $Y_1$ and the second eigenvalue of $Y_{d_0}$, we find that
	\[
	\ol{\theta} = (-\alpha \beta)^{d_0} \ol{\theta}, \]
	and so
	\[
	(-\alpha \beta)^{d_0} = 1. \]
\end{enumerate}	
\end{pf}


\subsection{} \label{sec3.21}
Next, we require the following elementary Lemma.  Recall that the \textbf{Hausdorff metric} $\textup{dist}_H(\cdot, \cdot)$ on the space $X = \{ A \subseteq \bbC:  \varnothing \neq A  \text{ compact}\}$ is defined by
\[
\textup{dist}_H(A, B) = \max (\max_{a \in A} \text{dist}(a, B), \max_{b \in B} \text{dist}(A, b)), \]
where, for $C \in X$, we set $\text{dist}(x, C) := \min \{ |x - c|: c \in C\}$.


\begin{lem}\label{lem3.22}
Let $4 \le d \in \bbN$ be an even integer and set $r = s = d_0 = \frac{d}{2}$.  Suppose that  $\alpha \ne \beta \in \bbT$ satisfy condition $\textbf{\emph{(P)}}$.
Define
\[
\begin{array}{rccc}
	f: & \bbT & \to & \bbR \\
		& \gamma & \mapsto & \underset{1\leq j\leq d_0}{\max}|(\alpha\beta)^{2j-1}-\gamma|.
\end{array} \]		
Then
\[\underset{\gamma\in \mathbb{T}}{\min}f(\gamma)=|e^{i\frac{d_0-1}{d_0}\pi}-1|.\]
\end{lem}

\begin{pf}
Set $\Gamma=\{(\alpha\beta)^{2j-1}:1\leq j\leq d_0\}$, and observe that for all $\gamma \in \bbT$, $f(\gamma) = \textup{dist}_H(\gamma, \Gamma)$.

Since $\alpha^{d_0} \beta^{d_0}\in \{-1,1\}$, we have $(\alpha\beta)^{2 d_0} = (\alpha \beta)^d =1$. Suppose that there exist $1\leq p<q\leq d_0$ such that
$(\alpha\beta)^{2p}=(\alpha\beta)^{2q}$.   Then $(\alpha\beta)^{2(q-p)}=1$. Hence $(\alpha\beta)^{q-p}\in \{-1,1\}$,
but $2(q-p)<2q\leq 2 d_0=d$, which contradicts condition \textbf{(P)}.
Thus,
\[(\alpha\beta)\cdot \Gamma=\{ ((\alpha\beta)^2)^j:1\leq j\leq d_0\}\]
has $d_0$ elements. Since $(\alpha\beta)^{2 d_0}=1$, all the elements listed in the above set are exactly the vertices of an inscribed regular $d_0$-polygon of $\mathbb{T}$.
Consequently, $\Gamma$ is the vertex set of an inscribed regular $d_0$-polygon of $\mathbb{T}$.

When $d_0$ is even, it follows from elementary plane geometry that $f(\gamma)$ attains its minimum value if and only if $\gamma$ is the middle point of an arc whose endpoints are two adjacent vertices of the inscribed regular $d_0$-polygon. Hence,
\[
f(\gamma) = \underset{\gamma\in \mathbb{T}}{\min}~\textup{dist}_H(\gamma, \Gamma) =|e^{i\frac{d_0-1}{d_0}\pi}-1|.\]

When $d_0$ is odd,  then, once again by elementary plane geometry, $f(\gamma)$ attains its minimum value if and only if $\gamma$ is a vertex of the inscribed regular $d_0$-polygon. Hence,
\[
f(\gamma) = \underset{\gamma\in \mathbb{T}}{\min}~\textup{dist}_H(\gamma, \Gamma)=|e^{i\frac{d_0-1}{d_0}\pi}-1|.\]
\end{pf}


\begin{thm} \label{thm3.23}
Let $4 \le d \in \bbN$ be an even integer and suppose that $d = r+s$ for some $s = r \in \bbN$.   Suppose that  $\alpha \ne \beta \in \bbT$ satisfy condition $\textbf{\emph{(P)}}$.
Let $d_0 := \frac{d}{2}$,  $V := \alpha I_r \oplus \beta I_s \in \cU(\bbC^d)$.
Then  $V \in \symthree(\bbC^d)$ and $\dim\, \cM_1 = 0 = \dim\, \cM_2$ in the standard decomposition of $J$ from Proposition~\ref{prop3.09} if and only if one of the following three conditions holds:
\begin{enumerate}
	\item[(a)]
	$d$ is divisible by $4$, $(\alpha \beta)^{d_0} = -1$, and
	\[
	|(\alpha \beta)^{2 j - 1} - 1|   < |\alpha - \beta|, \ \ \ 1 \le j \le d_0; \]
	\item[(b)]
	$d$ is not divisible by $4$,   and
	\[
	 |(\alpha \beta)^{2 j - 1} - 1| < |\alpha - \beta|, \ \ \ 1 \le j \le d_0; \text{ or } \]
	\item[(c)]
	$d$ is not divisible by $4$, $(-\alpha \beta)^{d_0} = 1$, and
	\[
	|e^{i\frac{d_0-1}{d_0}\pi}-1|<|\alpha-\beta|. \]
	\end{enumerate}

\end{thm}

\begin{pf}
	We begin by proving the necessity of these conditions.
As we have seen multiple times already, this implies that
	\[
	Z := V J = Y = \oplus_{j=1}^{d_0} Y_j, \]
	where $Y_j = \begin{bmatrix} \alpha m_j & \alpha n_j \\ \beta n_j & - \beta m_j \end{bmatrix}$ for some $m_j \in (-1, 1)$, $n_j \in (0, 1]$ for which $m_j^ 2 + n_j^2 =1$, $1 \le j \le d_0$.    In particular, $Z \in \symtwo(\bbC^d)$ and thus $Z \simeq Z^*$, implying that the non-real eigenvalues of $Z$ occur in conjugate pairs, counted according to multiplicity.

	
\bigskip

\noindent{\textbf{Case 1.}} $d \in 4 \bbN$.

If $d \in 4 \bbN$, then we cannot have $(- \alpha \beta)^{d_0} = 1$, otherwise $(- \alpha \beta)^{d/4} \in \{ -1, 1\}$, contradicting our minimality assumptions.

\bigskip

\begin{itemize}
	\item	{} {\textbf{Subcase 1.1}} Suppose that $d \in 4 \bbN$ and $\sigma(Z) \cap \{ -1, 1 \} \ne \varnothing$, say $\xi \in \sigma(Z) \cap \{ -1, 1\}$.

	Since $Y = Z$, we may apply Lemma~\ref{lem3.20} to conclude that
	\[
	\sigma(Z) = \sigma(Y) = \cup_{j=1}^{d_0} \Delta_j(\xi), \]
	and
	\[
	(-\alpha \beta)^{d_0} \in \{ -1, 1\}. \]
As noted above, if $d \in 4 \bbN$, then we cannot have $(-\alpha \beta)^{d_0} = 1$, and so $(\alpha \beta)^{d_0} = (-\alpha \beta)^{d_0} =-1$.

\bigskip

Note that for $1 \le j \le d_0$,
\begin{align*}
	|\textsc{Tr}(Y_j)|
		&= | (-\alpha \beta)^j \xi + \ol{(-\alpha \beta)^{j-1}} \xi| \\
		&= |(-\alpha \beta)^{2 j - 1}+ 1| \\
		&= | (\alpha - \beta) m_j| \\
		&< |\alpha - \beta|.
\end{align*}		

\bigskip

	\item{}{\textbf{Subcase 1.2}}. Now suppose that $d \in 4 \bbN$ and that $\sigma(Z) \cap \{ -1, 1\} = \varnothing$.

\bigskip

Let $\theta \in \sigma(Y)$.   By Lemma~\ref{lem3.20}, we conclude that -- after reindexing each of the $Y_j$'s if necessary -- we may assume that $\sigma(Y_j) = \Delta_j(\theta)$, $1 \le j \le d_0$.   Now the eigenvalues of $Y = Z$ appear in conjugate pairs, counted according to multiplicity, and the second eigenvalue of $Y_j$ is paired with the first eigenvalue of $Y_{j+1}$ for $1 \le j \le d_0-1$.  Since $\sigma(Y) \cap \bbR = \varnothing$, it follows that the second eigenvalue of $Y_{d_0}$ is paired with the first eigenvalue of $Y_1$, i.e.
\[
(-\alpha \beta)^{d_0} \ol{\theta} = \ol{\theta}. \]
But then $(- \alpha \beta)^{d_0} = 1$, and so
\[
(\alpha \beta)^{d_0/2} \in \{ -1, 1\}, \]
contradicting condition \textbf{(P)}.
	
\bigskip

\noindent{\textbf{Conclusion.}}
\begin{enumerate}
	\item[(a)]
	If $d \in 4 \bbN$, then $(-\alpha \beta)^{d_0} = -1$, and
\[
|(\alpha \beta)^{2 j - 1} - 1|=|(-\alpha \beta)^{2 j - 1} +1|  < |\alpha - \beta|, \ \ \ 1 \le j \le d_0. \]
\end{enumerate}
\end{itemize}

\bigskip

\noindent{\textbf{Case 2.}} $d \in 4 \bbN + 2$.

\begin{itemize}
	\item	{} {\textbf{Subcase 2.1}} $d \in 4 \bbN + 2$ and $\sigma(Z) \cap \{ -1, 1 \} \ne \varnothing$.

Now suppose that $d \in 4 \bbN + 2$, and suppose that $\xi \in \sigma(Z) \cap \{ -1, 1\} \ne \varnothing$.   By Lemma~\ref{lem3.20},
\[
\sigma(Z) = \sigma(Y) = \cup_{j=1}^{d_0} \Delta_j(\xi), \]
and
\[
(-\alpha \beta)^{d_0}  \in \{ -1,  1\}. \]

\bigskip

Note that for $1 \le j \le d_0$,
\begin{align*}
	|\textsc{Tr}(Y_j)|
		&= | (-\alpha \beta)^j \xi + \ol{(-\alpha \beta)^{j-1}} \xi| \\
		&= |(-\alpha \beta)^{2 j - 1}  +1| \\
		&= | (\alpha - \beta) m_j| \\
		&< |\alpha - \beta|.
\end{align*}

\bigskip

	\item{} {\textbf{Subcase 2.2}}. Now suppose that $d \in 4 \bbN +2$ and that $\sigma(Z) \cap \{ -1, 1\} = \varnothing$.

\bigskip

Let $\theta \in \sigma(Z) = \sigma(Y)$.   By Lemma~\ref{lem3.20},
\[
\sigma(Z) = \sigma(Y) = \cup_{j=1}^{d_0} \Delta_j(\theta), \]
and
\[
(-\alpha \beta)^{d_0}  = 1. \]
	
	As before, for $1 \le j \le d_0$,
\begin{align*}
	|\textsc{Tr}(Y_j)|
		&= | (-\alpha \beta)^j \ol{\theta} + \ol{(-\alpha \beta)^{j-1}} \theta| \\
		&= |(-\alpha \beta)^{2 j - 1} +\theta^2| \\
&= |(\alpha \beta)^{2 j - 1} -\theta^2|\\
		&= | (\alpha - \beta) m_j| \\
		&< |\alpha - \beta|.
\end{align*}		
By Lemma~\ref{lem3.22}, it follows that
\[
| e^{i \frac{d_0-1}{d_0} \pi} - 1| < | \alpha - \beta|. \]
	
\bigskip

\noindent{\textbf{Conclusion.}} If $d \in 4 \bbN + 2$, then $(-\alpha \beta)^{d_0} \in  \{ -1, 1\}$, and either
\begin{enumerate}
	\item[(b)]
	$\sigma(Y) \cap \{ -1, 1\} \ne \varnothing$, in which case
	\[
	 |(\alpha \beta)^{2 j - 1} - 1| = | (\alpha - \beta) m_j| < |\alpha - \beta|, \ \ \ 1 \le j \le d_0; \text{ or } \]
	\item[(c)]
	$\sigma(Y) \cap \{ -1, 1\} = \varnothing$, in which case 
	\[
	| e^{i \frac{d_0-1}{d_0} \pi} - 1 |   < |\alpha - \beta|, \ \ \ 1 \le j \le d_0. \]
\end{enumerate}

\end{itemize}


\vskip 0.5 cm

Next we will prove the sufficiency of these conditions. As always, this is essentially reverse-engineering the necessary conditions.

\begin{itemize}

\item{}  Suppose that (a) holds, namely:   $d \in 4 \bbN$, $(\alpha \beta)^{d_0} = -1$, and
\[
| (\alpha \beta)^{2j -1} - 1 | < | \alpha - \beta|, \ \ \ 1 \le j \le d_0. \]

Define
\[m_j=\frac{(-1)^j \ol{(\alpha \beta)}^{j-1}[(\alpha \beta)^{2j-1}-1]}{\alpha-\beta}=\frac{(-\ol{\alpha \beta})^{j-1} + (-\alpha \beta)^j }{\alpha-\beta},~~ 1 \le j \le  d_0.\]
Then
\[|m_j|=|\frac{(\alpha\beta)^{2j-1}-1}{\alpha-\beta}|<1, \ \ \ 1 \le j \le d_0.\]
By Lemma~\ref{lem3.07}, $m_j\in \bbR$,  whence $m_j\in (-1,1)$, and we define $n_j = \sqrt{1 - m_j^2}$, $1 \le j \le d_0.$
Now
\[
J:=\oplus_{j=1}^{d_0}\begin{bmatrix}  m_j &  n_j \\  n_j & - m_j \end{bmatrix}\]
is a symmetry with $\tr(J)=0$.
Consider
\[
V=\oplus_{j=1}^{d_0} \begin{bmatrix}  \alpha &  0 \\  0 & \beta \end{bmatrix} \simeq \alpha I_{d_0} \oplus \beta I_{d_0}.\]

For $1\leq j\leq d_0$, define
\[\lambda_1^{(j)}=(-\ol{\alpha \beta})^{j-1} , \lambda_2^{(j)}=(-\alpha \beta)^j .\]
Then it is straightforward to check that
\[\lambda_1^{(j)}+\lambda_2^{(j)}=(\alpha-\beta)m_j=\tr(\begin{bmatrix}  \alpha &  0 \\  0 & \beta \end{bmatrix}\begin{bmatrix}  m_j &  n_j \\  n_j & - m_j \end{bmatrix}),\]
and
\[\lambda_1^{(j)}\cdot\lambda_2^{(j)}=-\alpha\beta=\det(\begin{bmatrix}  \alpha &  0 \\  0 & \beta \end{bmatrix}\begin{bmatrix}  m_j &  n_j \\  n_j & - m_j \end{bmatrix}).\]
Hence, $\sigma(\begin{bmatrix}  \alpha &  0 \\  0 & \beta \end{bmatrix}\begin{bmatrix}  m_j &  n_j \\  n_j & - m_j \end{bmatrix})=\{\lambda_1^{(j)},\lambda_2^{(j)}\}$, $1\leq j\leq d_0$.

Now,
\[\lambda_2^{(j)}\lambda_1^{(j+1)}=1,~~\forall 1\leq j\leq d_0-1,\]
while
\[\lambda_1^{(1)}=1, \lambda_2^{(d_0)}= (-\alpha \beta)^{d_0}=-1.\]
Therefore, $VJ$ is a unitary whose non-real eigenvalues come in conjugate pairs.  In particular, $V J\in \symtwo(\bbC^d)$, and therefore  $V = (VJ) J \in \symthree(\bbC^d)$.


	\item{} Next suppose that (b) holds:  i.e., that $d \in 4 \bbN + 2$, and
	\[
	| (\alpha \beta)^{2j -1} - 1 | < | \alpha - \beta|, \ \ \ 1 \le j \le d_0. \]

	The proof of this case is essentially identical to that above.   The sole exception is that $\lambda_2^{(d_0)}$ can be either $1$ or $-1$,  depending upon the value of $(-\alpha \beta)^{d_0}$.


	\item{}  Finally, suppose that (c) holds:  i.e., that $d \in 4 \bbN + 2$, $(-\alpha \beta)^{d_0} = 1$, and
	\[
	|e^{i\frac{d_0-1}{d_0}\pi}-1|<|\alpha-\beta|. \]
	By Lemma~\ref{lem3.22}, there exists $\theta \in \mathbb{T}$,
	such that
	\[\underset{1\leq j\leq d_0}{\max}|(\alpha\beta)^{2j-1}-\theta^2|<|\alpha-\beta|.\]
	By continuity of the function $f$ of Lemma~\ref{lem3.22}, we may also
	assume that $\theta \notin \bbR$.
	For $1 \le j \le d_0$, set
	\[
	m_j :=\frac{(-1)^j\ol{(\alpha \beta)}^{j-1}\ol{\theta}[(\alpha \beta)^{2j-1}-\theta^2]}{\alpha-\beta}=\frac{(-\ol{\alpha \beta})^{j-1} \theta + (-\alpha \beta)^j \ol{\theta}}{\alpha-\beta},\]
	so that for each $1 \le j \le d_0$,
	\[
	|m_j|=|\frac{(\alpha\beta)^{2j-1}-\theta^2}{\alpha-\beta}|<1.\]
	By Lemma~\ref{lem3.07}, $m_j\in \bbR$,  whence $m_j\in (-1,1), 1 \le j \le  d_0$.
	Define $n_j=\sqrt{1-m_j^2}$, $1\leq j\leq d_0$. Then
	\[
	J:=\oplus_{j=1}^{d_0} \begin{bmatrix}  m_j &  n_j \\  n_j & - m_j \end{bmatrix}\]
	is a symmetry with $\tr(J)=0$.
	Consider
	\[
	V :=\oplus_{j=1}^{d_0}\begin{bmatrix}  \alpha &  0 \\  0 & \beta \end{bmatrix} \simeq \alpha I_{d_0} \oplus \beta I_{d_0}.\]
	Then $V J = \oplus_{j=1}^{d_0} \begin{bmatrix} \alpha m_j & \alpha n_j \\ \beta n_j & - \beta m_j \end{bmatrix}$.

	For $1\leq j\leq d_0$, define
	\[\lambda_1^{(j)}=(-\ol{\alpha \beta})^{j-1} \theta, \lambda_2^{(j)}=(-\alpha \beta)^j \ol{\theta}.\]
	As was the case earlier, it is straightforward to verify that for all $1 \le j \le d_0$,
	\[\lambda_1^{(j)}+\lambda_2^{(j)}=(\alpha-\beta)m_j=\tr(\begin{bmatrix}  \alpha m_j &  \alpha n_j \\  \beta n_j & - \beta m_j \end{bmatrix}),\]
	and
	\[\lambda_1^{(j)}\cdot\lambda_2^{(j)}=-\alpha\beta=\det(\begin{bmatrix}  \alpha m_j &  \alpha n_j \\  \beta n_j & - \beta m_j \end{bmatrix}).\]
	Hence, $\sigma(\begin{bmatrix}  \alpha m_j & \alpha n_j \\  \beta n_j & - \beta m_j \end{bmatrix})=\{\lambda_1^{(j)},\lambda_2^{(j)}\}$, $1\leq j\leq d_0$.

	Now,
	\[\lambda_2^{(j)}\lambda_1^{(j+1)}=1,~~\forall 1\leq j\leq d_0-1,\]
	while
	\[\lambda_1^{(1)}\cdot \lambda_2^{(d_0)}=\theta \cdot(-\alpha \beta)^{d_0} \ol{\theta}= (-\alpha \beta)^{d_0} =1.\]
	Therefore, $V J$ is a unitary whose eigenvalues are all non-real and come in conjugate pairs, implying that  $V J\in \symtwo(\bbC^d)$.   But then  $V\in \symthree(\bbC^d)$, completing the proof.
\end{itemize}
\end{pf}


\begin{thm}\label{thm3.24}
Let $4 \le d \in \bbN$ be an even integer and suppose that $d = r+s$ for some $r=s\in  \bbN$.   Suppose that  $\alpha \ne \beta \in \bbT$ satisfy condition $\textbf{\emph{(P)}}$.
Let $V := \alpha I_r \oplus \beta I_r \in \cU(\bbC^d)$.  Then  $V \in \symthree(\bbC^d)$ if and only if one of the following four conditions holds:
\begin{enumerate}
	\item[(a)]
	$d$ is divisible by $4$,  $(\alpha \beta)^{r} = -1$, and
	\[
	|(\alpha \beta)^{2 j - 1} - 1|   < |\alpha - \beta|, \ \ \ 1 \le j \le r; \]
	\item[(b)]
	$d$ is not divisible by $4$,  and
	\[
	 |(
\alpha \beta)^{2 j - 1} - 1| < |\alpha - \beta|, \ \ \ 1 \le j \le r;  \]
	\item[(c)]
	$d$ is not divisible by $4$,  $(\alpha \beta)^{r} = -1$, and
	\[
	|e^{i\frac{r-1}{r}\pi}-1|<|\alpha-\beta|. \]
\item[(d)] $d$ is not divisible by $4$, $(\alpha\beta)^{r} = -1$,
	\[
	| 1 - \alpha^{2 j+1}  \beta^{2 j - 1}|   < | \alpha - \beta|, \ \ \ \ \ 1 \le j \le r-1. \]	
	\end{enumerate}			
\end{thm}

\begin{pf}
Suppose that $V \in \symthree(\bbC^d)$.  Write $V = J_1 J_2 J$ where $J \in \symone(\bbC^d)$ is balanced, by Proposition~\ref{prop3.09} and Remark~\ref{rem3.10},
either $\dim\cM_1=0=\dim\cM_2$, or $\dim\cM_1=1=\dim\cM_2$.
Then by applying Theorem \ref{thm3.23} and Theorem \ref{thm3.17} respectively, we can deduce the desired conclusions.
\end{pf}


We now provide an alternate version of Theorem \ref{thm3.24}.

\begin{thm}\label{thm3.25}
Let $4 \le d \in \bbN$ be an even integer and suppose that $d = r+s$ for some $r=s\in  \bbN$.   Suppose that  $\alpha \ne \beta \in \bbT$ satisfy condition $\textbf{\emph{(P)}}$.
Let $V := \alpha I_r \oplus \beta I_r \in \cU(\bbC^d)$.  Then  $V \in \symthree(\bbC^d)$ if and only if one of the following three conditions holds:
\begin{enumerate}
	\item[(i)]
	$d$ is divisible by $4$,  $(\alpha \beta)^{r} = -1$, and
\[
	|e^{i\frac{r-1}{r}\pi}-1|<|\alpha-\beta|. \]
	\item[(ii)]
	$d$ is not divisible by $4$, and
	\[
	|e^{i\frac{r-1}{r}\pi}-1|<|\alpha-\beta|. \]
\item[(iii)] $d$ is not divisible by $4$, $(\alpha\beta)^{r} = -1$,
	\[
	| 1 - \alpha^{2 j+1}  \beta^{2 j - 1}|   < | \alpha - \beta|, \ \ \ \ \ 1 \le j \le r-1. \]	
	\end{enumerate}			
\end{thm}

\begin{pf}
We will prove this theorem by comparing  conditions (a), (b) (c), (d) of Theorem~\ref{thm3.24} and  conditions (i), (ii), (iii) herein.
\begin{itemize}
\item Suppose that $d$ is divisible by $4$,  $(\alpha \beta)^{r} = -1$. Then $r$ is even.
Recall the proof of Lemma \ref{lem3.22}, $\Gamma$ is the vertex set of an inscribed regular $r$-polygon, where $\Gamma=\{(\alpha\beta)^{2j-1}:1\leq j\leq r\}$. Since $1=(\alpha\beta)^{2r}$, $1$ is the  middle point of an arc whose endpoints are two adjacent vertices of this inscribed regular $r$-polygon. Hence \[|e^{i\frac{r-1}{r}\pi}-1|=\underset{1\leq j\leq r}{\max}|(\alpha\beta)^{2j-1}-1|.\]
Therefore, the condition (a) of Theorem \ref{thm3.24} is actually the same as condition (i).
\item Suppose condition (b) of Theorem \ref{thm3.24} holds. That is, $d$ is not divisible by $4$,  and
	\[
	 |(
\alpha \beta)^{2 j - 1} - 1| < |\alpha - \beta|, \ \ \ 1 \le j \le r.  \]

By Lemma \ref{lem3.22}, \[|e^{i\frac{r-1}{r}\pi}-1|\leq \underset{1\leq j\leq r}{\max}|(\alpha\beta)^{2j-1}-1|<|\alpha-\beta|.\]
Hence the condition (ii) holds. It is clear that the condition (ii) holds if the condition (c) of Theorem \ref{thm3.24} holds.
\item Suppose that condition (ii) holds. Then $r$ is odd.
\begin{enumerate}
\item If $(\alpha\beta)^r=1$. Then $\Gamma$ is the vertex set of an inscribed regular $r$-polygon, where $\Gamma=\{(\alpha\beta)^{2j-1}:1\leq j\leq r\}$. Since $1=(\alpha\beta)^r=(\alpha\beta)^{2\frac{r+1}{2}-1}\in \Gamma$, $1$ is a vertex of this inscribed regular $r$-polygon. Hence by the proof of Lemma \ref{lem3.22},\[|e^{i\frac{r-1}{r}\pi}-1|=\underset{1\leq j\leq r}{\max}|(\alpha\beta)^{2j-1}-1|<|\alpha-\beta|.\]
    Then the condition (b) of Theorem \ref{thm3.24} holds.
\item If $(\alpha\beta)^r=-1$. Then the condition (c) of Theorem \ref{thm3.24} holds.
\end{enumerate}
\item The condition (d) of Theorem \ref{thm3.24} is identical to the condition (iii).
\end{itemize}
\end{pf}


\begin{eg} \label{eg3.26}
Suppose that $d=2 d_0 =8$, $\varepsilon=\frac{\sqrt{2} \pi}{1000}$, and let
\[\alpha=e^{i (\pi+\varepsilon)},~~\beta=e^{i (\frac{\pi}{4}-\varepsilon)}.\] Then $(\alpha\beta)^4=-1$,
and it is routine to verify that if $i, j \in \bbN$ and  $\alpha^i \beta^j \in \{ -1, 1\}$, then $i+j \ge 8$. Let $V=\alpha I_4\oplus \beta I_4$, and note that $\det(V)=-1$.   We claim that $V\in \symthree(\bbC^8)$.

Observe that
\[|\sin \left(\frac{\arg{\alpha}-\arg{\beta}}{2}\right)|=|\sin (\frac{3\pi}{8}+\varepsilon)|=\sin (\frac{3\pi}{8}+\varepsilon)>\sin (\frac{3\pi}{8})=|\sin \left( \frac{\frac{4-1}{4}\pi}{2} \right )|.\]
By Lemma \ref{lem3.07},
\[|e^{i\frac{4-1}{4}\pi}-1|<|\alpha-\beta|.\]
Thus, according to Theorem \ref{thm3.25}(i), $V=\alpha I_4\oplus \beta I_4\in \symthree(\bbC^8)$.

We point out that the same argument works if we replace $\eps$ above by any  positive number $0 < \zeta < \frac{\sqrt{2} \pi}{1000}$, as long as $\frac{\zeta}{\pi}$ is irrational.
\end{eg}


\begin{eg} \label{eg3.27}
Suppose that $d=2 d_0 =6$, $\varepsilon = \frac{\sqrt{2} \pi}{1000}$, and let
\[\alpha=e^{i (\pi+\varepsilon)},~~\beta=e^{i (\frac{\pi}{3}-\varepsilon)}.\] Then $\alpha^3\beta^3=1$,
and as before, if $i, j \in \bbN$, then  $\alpha^i \beta^j \in \{ -1, 1\}$ implies that $i+j \ge 6$. Let $V=\alpha I_3\oplus \beta I_3$, then $\det(V)=1$.
 We claim that $V \in \symthree(\bbC^6)$.

To see this, note that
\[|\sin \left( \frac{\arg{\alpha}-\arg{\beta}}{2}\right) |=|\sin (\frac{\pi}{3}+\varepsilon)|=\sin (\frac{\pi}{3}+\varepsilon)>\sin (\frac{2\pi}{6})=|\sin \left(\frac{\frac{3-1}{3}\pi}{2}\right)|,\]
and thus by Lemma \ref{lem3.07},
\[|e^{i\frac{3-1}{3}\pi}-1|<|\alpha-\beta|.\]

It follows from Theorem \ref{thm3.25}(ii) that $V=\alpha I_3\oplus \beta I_3\in \symthree(\bbC^6)$.
\end{eg}


\begin{eg} \label{eg3.28}
Assume that $d=2 r =2s=6$, and set $\eps = \frac{\sqrt{2} \pi}{1000}$.   Set
\[
\alpha=e^{i (\pi+\frac{\pi}{6}+\varepsilon)},~~\beta=e^{i (\frac{\pi}{2}-\varepsilon)}.\]
Then $\alpha^3\beta^3=-1$,  $\alpha^i \beta^j \in \{ -1, 1\}$ implies that $i+j \ge 6$. Let $V=\alpha I_3\oplus \beta I_3$, then $\det(V)=-1$. We claim that $V\in \symthree(\bbC^6)$.
Note that
\[|\sin \left( \frac{\arg{\alpha}-\arg{\beta}}{2} \right)|=|\sin (\frac{\pi}{3}+\varepsilon)|=\sin (\frac{\pi}{3}+\varepsilon),\]
\[|\sin \left( \frac{\arg(\alpha^3\beta)}{2} \right)|=|\sin \left( \frac{3(\pi+\frac{\pi}{6}+\varepsilon)+(\frac{\pi}{2}-\varepsilon)}{2} \right)|=\sin (\varepsilon),\]
and
\[|\sin \left( \frac{\arg(\alpha^5\beta^3)}{2} \right)|=|\sin \left( \frac{5(\pi+\frac{\pi}{6}+\varepsilon)+3(\frac{\pi}{2}-\varepsilon)}{2}\right)|=\sin (\frac{\pi}{3}-\varepsilon).\]
Therefore, by Lemma \ref{lem3.07},
\[\frac{|\alpha^3\beta-1|}{|\alpha-\beta|}=\frac{|\sin \left(\frac{\arg(\alpha^3\beta)-0}{2} \right)|}{|\sin \left( \frac{\arg{\alpha}-\arg{\beta}}{2}\right)|}<1,\]
and
\[\frac{|\alpha^5\beta^3-1|}{|\alpha-\beta|}=\frac{|\sin \left( \frac{\arg(\alpha^5\beta^3)-0}{2}\right)|}{|\sin \left(\frac{\arg{\alpha}-\arg{\beta}}{2} \right)|}<1.\]
So by Theorem~\ref{thm3.25}(iii), $V=\alpha I_3\oplus \beta I_3\in \symthree(\bbC^6)$.

Alternatively, note that
\[|\sin \left( \frac{\arg{\alpha}-\arg{\beta}}{2} \right)|=\sin (\frac{\pi}{3}+\varepsilon)>\sin (\frac{2\pi}{6})=|\sin\frac{\frac{3-1}{3}\pi}{2}|.\]
By Lemma~\ref{lem3.07},
\[|e^{i\frac{3-1}{3}\pi}-1|<|\alpha-\beta|.\]

By  Theorem~\ref{thm3.25}(ii), $V=\alpha I_3\oplus \beta I_3\in \symthree(\bbC^6)$.
\end{eg}


\bigskip

It is reasonable to ask which compact sets can occur as the spectra of elements of $\symthree(\hilb)$, whether $\hilb$ is finite- or infinite-dimensional. The results of this section provide an interesting set of examples.

\bigskip

\begin{lem} \label{lem3.29}
Let $p\geq 3$ be a prime number, and $\alpha$ be a primitive $p$-th root of unity.   Set $\beta=-\alpha$, $r=\frac{p+1}{2}$, and $s=\frac{p-1}{2}$. Then $V=\alpha I_r\oplus \beta I_s\in \symthree(\bbC^{p})$.
\end{lem}

\begin{pf}
Note that $\alpha^r\beta^s=(-1)^s\alpha^{p}=(-1)^s$. Suppose that $i,j\in \bbN$ and that $\alpha^i\beta^{j}\in \{-1,1\}$. Then $\alpha^{i+j}\in \{-1,1\}$.
Since $p\geq 3$ is a prime number, and $\alpha$ is a primitive $p$-th root of unity, it follows that $i+j\geq p$. That is,  $\alpha \ne \beta \in \bbT$ satisfy condition \textup{\textbf{(P)}}.
In particular, 	
\[| 1 - (\alpha \beta)^{2 j - 1}|   < 2=|\alpha-\beta|, 1 \le j \le s. \]	
By Theorem \ref{thm3.12}, $V=\alpha I_r\oplus \beta I_s\in \symthree(\bbC^{p})$.
\end{pf}


\begin{prop} \label{prop3.30}
Let $\cH$ be a separable, infinite-dimensional complex Hilbert space and $E :=\{\alpha\in \bbC: \alpha\neq 1,~~\alpha^p=1$ for some prime number $p\geq 3\}$.
Suppose that
$K$ is a non-empty compact subset of $\bbT$ which is symmetric with respect to the origin $($i.e. $K=-K$$)$ and that $K\cap E$ is dense in $K$.
Then there exists $V\in \symthree(\cH)$ such that $\sigma(V)=K$.
\end{prop}

\begin{pf}
Since $\{-1,1\}\cap E=\varnothing$, it follows that $K\neq \{-1,1\}$.
Write $K\cap E = \{\alpha_k\}_{k=1}^\infty$ (where we allow $\alpha_i = \alpha_j$ with $i \ne j$).

For $k\in \bbN$, there exists a prime number $p_k\geq 3$, such that $\alpha^{p_k}=1$.  Set $r_k=\frac{p_k+1}{2}$ and $s_k=\frac{p_k-1}{2}$.
By Lemma \ref{lem3.29}, $V_k=(\alpha_k)I_{r_k}\oplus (-\alpha_k)I_{s_k}\in \symthree(\bbC^{p_k})$, $k\in \bbN$.

The operator $V := \oplus_{k =1}^\infty V_k \in \bofh$ is unitary, and it is routine to verify that $V \in \symthree(\hilb)$ because each $V_k \in \symthree(\bbC^{p_k})$.
Since $K=-K$,
\[K=\overline{K\cap E}=\overline{\{\alpha_k:k\geq \bbN\}}\subset \sigma(V)=\overline{\{\alpha_k,-\alpha_k:k\geq \bbN\}}\subset \overline{K\cup (-K)}=K.\]
Therefore, $\sigma(V)=K$.
\end{pf}



\vskip 1 cm
\section{General unitary matrices in $\symthree(\bbC^d)$} \label{section4}


\subsection{} \label{sec4.01}
The previous section illustrates just how complicated the structure of $\symthree(\bbC^d)$ is.   At the moment, a complete description of those unitary matrices which lie in this set seems beyond reach.   Despite this, a number of partial results have been obtained.   For example, the second author~\cite{Radjavi1969c} has shown that if $\Psi_m := \{ e^{2 \pi i t}  : \frac{m-1}{4} < t < \frac{m}{4} \}$, $1 \le m \le 4$, then no unitary operator $V \in \bofh$ whose spectrum is contained in a single $\Psi_m$ can lie in $\symthree(\hilb)$.

\smallskip

In this section, we obtain a number of results which impose necessary conditions for a unitary operator $V$ to lie in $\symthree(\bbC^d)$.   Our first conditions require the notion of the \emph{length} of a product of reflections.


\begin{defn} \label{def4.02}
By a  \textbf{reflection} in $\cU(\bbC^d)$ we mean a symmetry which is a rank-one perturbation of the identity operator $I_d$.
We denote by $\textsc{ref}(\bbC^d)$ the set of all reflections in $\cU(\bbC^d)$, and we define
\[
\prod(\textsc{ref}(\bbC^d)) := \{ R_1 R_2 \cdots R_m :m \ge 1 \text{ and } R_j \in \textsc{ref}(\bbC^d), 1 \le j \le m \}. \]
Given $W \in \prod(\textsc{ref}(\bbC^d))$, we define the \textbf{length} of $W$ to be
\[
\ell({W}) := \min\{ m: W = R_1 R_2 \cdots R_m:  R_j \in \textsc{ref}(\bbC^d), 1 \le j \le m \}. \]
\end{defn}


\subsection{} \label{sec4.03}
In~\cite[Theorem~9]{Radjavi1975}, the second author showed that $l(W)\leq 2d-1$ for any unitary matrix $W \in \cU(\bbC^d)$ whose determinant is real.

Let $W \in \cU(\bbC^d)$.   Then there exist uniquely determined real numbers
\[
0 \leq \alpha_1 \leq \alpha_2 \le  \cdots \leq \alpha_d < 2\pi\]
such that $W$ is unitarily equivalent to the diagonal matrix $D := \textsc{diag}(e^{i \alpha_1}, e^{i \alpha_2}, \ldots, e^{i \alpha_d})$.
We refer to $\alpha_1, \alpha_2, \ldots,\alpha_d$ as the \textbf{angles} of $W$, and we define
\[
\kappa (W)=\frac{\sum_{i=1}^d \alpha_i}{\pi}.\]


The following theorem is due to Djokovi\'{c} and Malzan \cite{DjokovicMalzan1979}.

\begin{thm} \emph{\textbf{[Djokovi\'{c} and Malzan]}} \label{thm4.04}
Let $d \in \bbN$.  If $W \in  \cU(\bbC^d)$ and $\det(W)\in \{-1,1\}$, then
\[
\ell(W)=\max\{ \kappa (W), \kappa (W^*)\}. \]
\end{thm}


\begin{lem}\label{lem4.05}
If $W \in \symtwo(\bbC^d)$, then $\ell(W)\leq d$.
\end{lem}

\begin{pf}
Since $W\in \symtwo(\bbC^d)$, $W$ is unitarily equivalent to $W^*$. As the length $\ell(W)$ of $W$ is invariant under unitary equivalence, we may assume that
\[
W=\begin{bmatrix} e^{i\alpha_1}&&&&&\\ &e^{i(2\pi-\alpha_1)}&&&&\\&& \ddots &&&\\ &&& \ddots && \\ &&&& e^{i\alpha_r}&\\ &&&&& e^{i(2\pi-\alpha_r)}\\\end{bmatrix}\oplus I_s\oplus (-I_t),\]
where $0<\alpha_1\leq \cdots \leq \alpha_r<\pi$, and $2r+s+t=d$.
Note that $\kappa(W)=\kappa(W^*)$. Hence, by Theorem~\ref{thm4.04},
\[
\ell(W)=\kappa(W)=\frac{r\cdot 2\pi+t\cdot \pi}{\pi}\leq 2r+t\leq d.\]
\end{pf}


\begin{prop}\label{prop4.06}
If $d \in \bbN$ and $W\in \textsc{Sym}_3(\bbC^d)$, then $\ell(W)\leq d +\lfloor \frac{d}{2} \rfloor$.
\end{prop}

\begin{pf}
Write $W = J_1 J_2 J_3$, where $J_j \in \symone(\bbC^d)$, $j=1,2,3$.
Since $W =J_1(-J_2)(-J_3)$ and
\[
\min\{\rank(J_3-I_d), \rank(-J_3-I_d)\}\leq \lfloor \frac{d}{2} \rfloor,\]
we may assume without loss of generality that $\rank(J_3-I_d)\leq \lfloor \frac{d}{2} \rfloor$. Thus, $\ell(J_3)\leq \lfloor \frac{d}{2} \rfloor$.
According to the definition of the length of $W$ and Lemma \ref{lem4.05},
\[
\ell(W)\leq \ell(J_1 J_2)+ \ell(J_3)\leq d+\lfloor \frac{d}{2} \rfloor.\]
\end{pf}


The above Proposition allows us to provide an alternate proof to a finite-dimensional version of  \cite[Corollary 2]{Radjavi1969c}.  For $1 \le m \le 4$, let
\[
\Psi_m := \left \{ e^{i \theta}: \frac{(m-1) \pi}{2} < \theta < \frac{m \pi}{2} \right \}. \]

\begin{cor} \label{cor4.07}
Let $d \in \bbN$.   No unitary matrix whose spectrum is a subset of one of the four open arcs $\Psi_m$, $1 \le m \le 4$  lies in $\symthree(\bbC^d)$.
\end{cor}

\begin{pf}
Suppose that $W \in\cU(\bbC^{d})$.  Then it routine to see that the following conditions are equivalent;
\begin{center}
\begin{tabular}{lcl} 
     $\bullet  \ W \in \symthree(\bbC^d)$  & \ \ \ \ \  \ \ \ \ \ & $\bullet  \   -W \in \symthree(\bbC^d)$\\
     $\bullet  \ W^* \in \symthree(\bbC^d)$ & \ \ \ \ \ \ \ \ \ \  &  $\bullet  \ -W^* \in \symthree(\bbC^d)$ \\
\end{tabular}
\end{center}

If $\det(W)\notin \{-1,1\}$, then $W$ is not a finite product of symmetries, and so we may assume that $\det(W)\in \{-1,1\}$.
\begin{enumerate}
	\item[(i)]
	If $\sigma(W)\subseteq \Psi_4$, then by Theorem~\ref{thm4.04},
	\[
	\ell(W)\geq \kappa(W) > \frac{\frac{3\pi}{2}\cdot d}{\pi}= \frac{3d}{2}\geq d+\lfloor \frac{d}{2} \rfloor.\]
	By Proposition~\ref{prop4.06}, $W \notin \symthree(\bbC^d)$.
	\item[(ii)]
	If $\sigma(W)\subseteq  \Psi_1$, then $\sigma(W^*)\subseteq \Psi_4$.   By (i),  $W^*\notin \symthree(\bbC^d)$, whence $W\notin \symthree(\bbC^d)$.
	\item[(iii)]
	If $\sigma(W)\subseteq  \Psi_2$, then $\sigma(-W)\subseteq \Psi_4$.   By (i), $-W \notin  \symthree(\bbC^d)$, whence $W\notin \symthree(\bbC^d)$.
	\item
	If $\sigma(W)\subseteq  \Psi_3$, then $\sigma(-W^*)\subseteq \Gamma_4$.  By (i), $-W^* \notin  \symthree(\bbC^d)$, whence $W\notin \symthree(\bbC^d)$.
\end{enumerate}
\end{pf}


Proposition~\ref{prop4.06} is a useful device for producing examples of unitary operators in $\cU(\bbC^d)$ whose determinant is real, but which \emph{do not} lie in $\symthree(\bbC^d)$.


\begin{prop}\label{prop4.08}
Let $3 \le d \in \bbN$, $W\in  \cU(\bbC^d)$ and suppose that $\det(W)\in \{-1,1\}$.  Suppose furthermore that $W\simeq \oplus_{j=1}^k (e^{i\theta_j}I_{m_j})$,
where $\pi<\theta_1<\cdots<\theta_k<2\pi$. Then
\[
\ell(W)+\ell(-W)=3d.\]
In particular, if $d$ is odd or $\ell(W) \ne d +\lfloor \frac{d}{2} \rfloor$, then $W \notin \symthree(\bbC^d)$.
\end{prop}

\begin{pf}
First, observe that
\[
W^*\simeq \oplus_{j=1}^k (e^{i(2\pi-\theta_j)}I_{m_j}),\ \ \ \ \ -W \simeq \oplus_{j=1}^k (e^{i(\theta_j-\pi)}I_{m_j}), \ \text{ and } \ \ \ -W^*\simeq \oplus_{j=1}^k (e^{i(3\pi-\theta_j)}I_{m_j}).\]
Thus
\begin{align*}
\ell(W)
	&=\max\{ \kappa(W), \kappa(W^*)\} \\
	&=\max\{\frac{\sum_{j=1}^k m_j\theta_j}{\pi},\frac{\sum_{j=1}^k m_j(2\pi-\theta_j)}{\pi}\} \\
	&=\frac{\sum_{j=1}^k m_j\theta_j}{\pi},
\end{align*}	
while
\begin{align*}
\ell(-W)
	&=\max\{ \kappa (-W), \kappa(-W^*)\}  \\
	&=\max\{\frac{\sum_{j=1}^k m_j(\theta_j-\pi)}{\pi},\frac{\sum_{j=1}^k m_j(3\pi-\theta_j)}{\pi}\}  \\
	&=\frac{\sum_{j=1}^k m_j(3\pi-\theta_j)}{\pi}.
\end{align*}	
Thus
\[
\ell(W)+\ell(-W)=\frac{\sum_{j=1}^k m_j\theta_j+\sum_{j=1}^k m_j(3\pi-\theta_j)}{\pi}=\frac{3\pi(\sum_{i=1}^k m_j)}{\pi}=3d.\]

\begin{itemize}
	\item{}
	If $d$ is odd, then $2(d+\lfloor \frac{d}{2} \rfloor)<3d$, so one of $\ell(W)$ and $\ell(-W)$ is larger than $d+\lfloor \frac{d}{2} \rfloor$.   By Proposition~\ref{prop4.06}, $W \notin \symthree(\bbC^d)$;
	\item{}
	if $d$ is even and $\ell(W)\neq d+\lfloor \frac{d}{2} \rfloor$, then either $\ell(W)$ or $\ell(-W)$ is  larger than $d+\lfloor \frac{d}{2} \rfloor$.   By Proposition~\ref{prop4.06} once again, $W \notin \symthree(\bbC^d)$.
\end{itemize}
This completes the proof.	
\end{pf}


\begin{rem} \label{rem4.09}
Suppose that $d \in \bbN$ is even and write $d = 2 d_0$.  Suppose that $W \in \cU(\bbC^d)$ is of the form
\[
W \simeq \ttt{diag}(e^{i \theta_1}, e^{i \theta_2}, \ldots, e^{i \theta_d}), \text{ where } \pi < \theta_1 \le \theta_2 \le \cdots \le \theta_d < 2 \pi. \]
If $W \in \symthree(\bbC^d)$, then  $\mathrm{det}(W) = (-1)^{d_0}$.

\bigskip

Indeed, suppose first that $d_0$ is odd, say $d_0 = 2 k + 1$ for some $k \in \bbN$.   Then $d = 4k + 2$.   By Proposition~\ref{prop4.08}, if $W \in \symthree(\bbC^d)$, then

\[
\ell(W) = \ell(-W) = d + \lfloor\frac{d}{2}\rfloor = \frac{3d}{2} = 3 (2k +1). \]
But then
\[
\mathrm{det}(W) = \prod_{j=1}^d e^{i \theta_j} = e^{i \sum_{j=1}^d \theta_j} = e^{\frac{3  d \pi i}{2}} = e^{3 (2k +1) \pi i} = -1 = (-1)^{d_0}. \]

\bigskip

If $d_0$ is even, say $d_0 = 2k$ for some $k \in \bbN$, then $d = 4k$ and again,
\[
\ell(W) = \ell(-W) = d + \lfloor\frac{d}{2}\rfloor = \frac{3d}{2} = 6k. \]
In this instance,
\[
\mathrm{det}(W) = e^{\frac{3 d \pi i}{2}} = e^{6 k \pi i} = 1 = (-1)^{d_0}. \]

\end{rem}

\begin{eg} \label{eg4.10}
Let $\eps := \frac{1}{100}$, and set $\alpha_1 := e^{\frac{17}{12} \pi i}$, $\alpha_2 := e^{( \frac{17}{12}\pi  + \eps)i}$, $\alpha_3 := e^{( \frac{17}{12} \pi - \eps)i}$.   Finally, let $\beta :=  e^{\frac{15}{8} \pi i}$, and define
\[
V = \ttt{diag} (\alpha_1, \alpha_2, \alpha_3, \beta, \beta) \in \cU(\bbC^5), \]
and let $W := V \oplus V \in \cU(\bbC^{10})$.
\bigskip

First observe that the eigenvalues of $V$  all lie in the lower half-plane, and so the analysis of the above Remark applies.  Next,
 note that $\mathrm{det}(V) = 1$.    Thus  $\mathrm{det}(W) = 1$ as well.   Observe furthermore that $10/2$ is odd.   By the above Remark with $d=10$, $d_0=5$, if $W \in \symthree(\bbC^{10})$, then $\mathrm{det}(W) = (-1)^5 = -1$.   Hence $W \not \in \symthree(\bbC^{10})$.
 \end{eg}


\bigskip

We finish this section with a minor result concerning operators in $\symtwo(\bbC^3)$ with up to three distinct eigenvalues.

\smallskip

\begin{prop} \label{prop5.20}
Suppose that $\alpha, \beta, \gamma \in \bbT$, $\alpha \, \beta \, \gamma \in \{ -1, 1\}$, and that $\varrho := \alpha + \beta + \gamma \in \bbR$.
Then $\{ \alpha, \beta, \gamma \} \cap \{ -1, 1 \} \ne \varnothing$, and
\[
\ttt{diag} (\alpha, \beta, \gamma) \in \symtwo(\bbC^3). \]
\end{prop}

\begin{pf}
Assume, without loss of generality, that $\alpha \, \beta\, \gamma = 1$. (Otherwise, we multiply everything by $-1$.) Note that then

\[
\alpha \beta + \beta \gamma + \gamma \alpha = \ol{\gamma} + \ol{\alpha} + \ol{\beta} = \ol{\alpha + \beta + \gamma} = \ol{\varrho} = \varrho. \]

Thus $p(x) := (x - \alpha) (x-\beta) (x-\gamma)$ is of the form
\[
p(x) = x^3  - \varrho x^2 + \varrho x - 1 = 0,\]
and so $1$ is a root of $p(x)$.

By relabelling if necessary, we may assume that $\alpha = 1$.   Since $\alpha \, \beta \, \gamma = 1$, it then follows that $\gamma = \ol{\beta}$.

But then $U := \ttt{diag} (\alpha, \beta, \gamma)$ is unitarily equivalent to $U^*$, and so $U \in \symtwo(\bbC^3)$.
\end{pf}



\vskip 0.5 cm

\section{Symmetric orbits} \label{section5}


\subsection{}  \label{sec5.01}
One of the  most important notions of equivalence of operators on a finite-dimensional space is that of unitary equivalence.    We have already seen that every unitary operator $V \in \cU(\bbC^d)$ with $\det(V)=\pm 1$ can be written as a product of four symmetries;  i.e. $\cU(\bbC^d) = \symfour(\bbC^d)$, and this leads us to formulate the concept of the $\textsc{Sym}_k$-orbit of $T$, $k =1, 2, 3, 4$, defined as follows.


\begin{defn} \label{defn5.02}
Given $d \in \bbN$, $T \in \cB(\bbC^d)$ and $1 \le k \le 4$,   we define
\[
\textsc{Sym}_k(T) := \{ W^* T W: W \in \textsc{Sym}_k(\bbC^d) \}, \]
and refer to this as the $\textsc{Sym}_k$-\textbf{orbit} of $T$.
\end{defn}

We note that when $k =4$, $\mathbb{T}\cdot\symfour(\bbC^d) = \cU(\bbC^d)$, and thus $\symfour(T) = \cU(T)$ for all $T \in \cB(\bbC^d)$.


\begin{rem} \label{rem5.03}
Fixing $1 \le k \le 3$, consider the relation $T_1 \equiv_k T_2$ if $T_2 \in \textsc{Sym}_k(T_1)$.   It is not hard to see that this relation is both reflexive and symmetric.
Because $\textsc{Sym}_k(\bbC^d)$ is not a group, however, there is no reason to believe that it should be transitive.
\end{rem}


\subsection{} \label{sec5.04}
Our main goal in this section is to prove that if $d \in \bbN$ and $T \in \cB(\bbC^d)$, then $\symone(T) = \cU(T)$ if and only if $T$ is normal and has at most two eigenvalues.

As it turns out, determining when $\symtwo(T) = \cU(T)$ or $\symthree(T)= \cU(T)$ appears to be quite difficult, even if $T$ is a normal (and therefore diagonalisable) operator.   We shall have more to say about this below.

\smallskip

Before beginning our analysis of the $\textsc{Sym}_k$-orbit of an operator, we need to define the notion of the \emph{unitary commutant} of a subset of $\bofh$.    At this stage, the definitions and arguments apply to operators acting upon both finite-dimensional and infinite-dimensional Hilbert spaces, and we formulate our first results in full generality.


\begin{defn} \label{definition5.05}
Let $\hilb$ be a complex Hilbert space and $\varnothing \ne \cS \subseteq \bofh$.   We define the \textbf{unitary commutant} of $\cS$ to be
\[
\uncom(\cS) := \{ U \in \cU(\hilb) : U T = T U \text{ for all } T \in \cS \} = \cU(\hilb) \cap \cS^\prime. \]
When $\cS = \{ T \}$ is a singleton set, we abbreviate this to $\uncom(T)$.
\end{defn}

It is routine to verify that $\uncom(\cS)$ is a group.


\smallskip

The following Proposition is rather elementary, but it will prove very useful.

\begin{prop} \label{prop5.06}
Let $\hilb$ be a complex Hilbert space, $1 \le k \le 4$, and suppose that $T \in \bofh$.  If $R \in \symk(\hilb)$, $U \in \cU(\hilb)$ and $R^* T R = U^* T U$, then $R^* X R = U^* X U$ for all $X \in C^*(T)$.

As a consequence, If $\symk(T) = \cU(T)$, then $\symk(X) = \cU(X)$ for all $X \in C^*(T)$.
\end{prop}

\begin{pf}
Let $U \in \cU(\hilb)$ and choose $R \in \symk(\hilb)$ such that $R^* T R = U^* T U$.    Then $R^* T^* R = (R^* T R)^* = (U^* T U)^* = U^* T^* U$.   A routine calculation  shows that $p(x,y)$ is a polynomial in two non-commuting variables, then
\[
R^* p(T, T^*) R = U^* p(T, T^*) U.\]

In the finite-dimensional setting, every element of $C^*(T)$ is of the form $p(T, T^*)$ for some such polynomial, completing the proof.

In the infinite-dimensional setting,  $X \in C^*(T)$ implies that $X = \lim_n p_n(T, T^*)$ for some sequence $p_n(x,y)$ of such polynomials, and thus
\begin{align*}
R^* X R = R^*(\lim_n p_n(T, T^*)) R
	&= \lim_n R^* (p_n(T, T^*) R \\
	&= \lim_n U^* (p_n(T, T^*)) U \\
	&= U^* (\lim_n p_n(T, T^*)) U \\
	&= U^* X U.
\end{align*}
The last statement clearly follows from this.
\end{pf}


\begin{lem} \label{lem5.07}
Let $\hilb$ be a complex Hilbert space, $T \in \bofh$, and $1 \le k \le 4$.  The following statements are equivalent.
\begin{enumerate}
	\item[(a)]
	$\cU(\hilb) = \uncom(T) \cdot \symk(\hilb) \cdot \uncom(T) := \{ G J H: G, H \in \uncom(T), J \in \symk(\hilb)\}$.	
	\item[(b)]
	$\cU(\hilb) = \uncom(T) \cdot \symk(\hilb) : = \{ G J : G \in \uncom(T) , J \in \symk(\hilb)\}$.
	\item[(c)]
	$\cU(T) = \symk(T)$.
\end{enumerate}	
\end{lem}

\begin{pf}
\begin{enumerate}
	\item[(a)] implies (b). \ \ \
	Suppose that $U = G J H$, where $G, H \in \uncom(T)$ and $J \in \symk(\hilb)$.    Then $U = (G H) (H^* J H)$, and $GH \in \uncom(T)$, $H^* J H \in \symk(T)$, so that $U \in  \uncom(T) \cdot \symk(\hilb)$.   Thus $\cU(\hilb) \subseteq  \uncom(T) \cdot \symk(\hilb)$, and the reverse inclusion is trivial.
	\item[(b)] implies (c). \ \ \
	Suppose that $\cU(\hilb) = \uncom(T) \cdot \symk(\hilb)$.   Given $U \in \cU(\hilb)$, choose $V \in \uncom(T)$ and $J \in \symk(\hilb)$ such that $U = V J$.
	
	 Then
	\[
	U^* T U = J^* V^* T V J = J^* V^* V T J = J^* T J \in \symk(T). \]
	Since $U \in \cU(\hilb)$ was arbitrary, $\cU(T) \subseteq \symk(T)$, whence $\cU(T) = \symk(T)$.
	\item[(c)] implies (a). \ \ \
	Now suppose that $\cU(T) = \symk(T)$.   Given $U \in \cU(\hilb)$, find $J_U \in \symk(\hilb)$ such that
	\[
	U^* TU = J_U^* T J_U. \]
	Then $T (U J_U^*) = (U J_U^*) T$,  so that $G := U J_U^*  \in \uncom(T)$.  From this we see that
	$U = G J_U I \in \uncom(T) \cdot \symk(\hilb) \cdot \uncom(T)$.   Since $U \in \cU(\hilb)$ was arbitrary, we are done.
\end{enumerate}	
\end{pf}




\begin{lem} \label{lem5.08}
Let $\hilb$ be a complex Hilbert space, $m\ge 3$, and suppose that $\hilb = \oplus_{j=1}^m \hilb_j$, where $\hilb_j \ne \{ 0\}$, $1 \le j \le m$.   Denote by $P_j$ the orthogonal projection of $\hilb$ onto $\hilb_j$ for all $j$.   Then there exists $U \in \cU(\hilb)$ such that
\[
\mathrm{rank}\, P_3 U P_1 \ne \mathrm{rank}\, P_1 U P_3. \]
\end{lem}

\begin{pf}
Choose unit vectors $e \in \hilb_1$, $f \in \hilb_2$ and $g \in \hilb_3$ and set $\cM := \mathrm{span} \{ e, f, g\}$.  Define the unitary operator $U = V \oplus W$, where $V = g \otimes e^* + e \otimes f^* + f \otimes g^* \in \cB(\cM)$, and $W = I|_{\cM^\perp}$.   It is routine to verify that $U$ is unitary.

Furthermore, $P_3 U P_1 = g \otimes e^*$ has rank one, while $P_1 U P_3 = 0$.
\end{pf}


\begin{prop} \label{prop5.09}
Let  $\hilb$ be a complex Hilbert space and $T \in \cB(\hilb)$.   If there exists a normal operator $N \in C^*(T)$ such that $\sigma(N)$ has at least three components, then
\[
\cU(T) \ne \symone(T). \]
\end{prop}

\begin{pf}
We shall argue by contradiction.   Suppose to the contrary that $N \in C^*(T)$, $\sigma(N)$ has at least three components $\Omega_1, \Omega_2, \Omega_3$, and that $\cU(T) = \symone(T)$.  Write $\hilb = \hilb_1 \oplus \hilb_2 \oplus \hilb_3$, where $\hilb_1 := \hilb(N; \Omega_1)$, $\hilb_2 := \hilb(N; \Omega_2)$, and $\hilb_3 := (\hilb_1 \oplus \hilb_2)^\perp$ are the spectral subspaces for $N$ arising from the Spectral Theorem.   With respect to the decomposition $\hilb = \hilb_1 \oplus \hilb_2 \oplus \hilb_3$, we may write
\[
N = N_1 \oplus N_2 \oplus N_3, \]
where $\sigma(N_i) \cap \sigma(N_j) = \varnothing$, $1 \le i \ne  j \le 3$.

By Lemma~\ref{lem5.08}, there exists $U \in \cU(\hilb)$ such that with
\[
\mathrm{rank}\, P_3 U P_1 \ne \mathrm{rank}\, P_1 U P_3, \]
where $P_j$ denotes the orthogonal projection of $\hilb$ onto $\hilb_j$, $1 \le j \le 3$.

Now, by Lemma~\ref{lem5.07},
\[
\cU(\hilb) = \uncom(T) \cdot \symone(T). \]
Choose $V \in \uncom(T)$ and $J \in \symone(\hilb)$ such that $U = V J$.   By Fuglede's Theorem, since $V T = T V$, we find that $V \in C^*(T)^\prime$, and hence $V N = N V$.  But the fact that the spectra of $N_1, N_2$ and $N_3$ are disjoint then implies that  $V = V_1 \oplus V_2 \oplus V_3$, where each $V_j \in \cU(\hilb_j)$, $1 \le j \le 3$.   Write $J= [J_{ij}]$ relative to the same decomposition of $\hilb$, and observe that the fact that $J=J^*$ implies that $\mathrm{rank}\, J_{13} = \mathrm{rank}\, J_{13}^* = \mathrm{rank}\, J_{31}$.

Thus
\[
\mathrm{rank}\, P_3 (V J) P_1 = \mathrm{rank}\, V_3 J_{31} = \mathrm{rank}\, J_{31} = \mathrm{rank}\, J_{13} = \mathrm{rank}\, V_1 J_{13} = \mathrm{rank}\, P_1 (VJ) P_3, \]
from which we conclude that $V J \ne U$, a contradiction.
It follows that $\cU(T) \ne \symone(T)$.
\end{pf}


Recall that if $d \in \bbN$ and $T \in \cB(\bbC^d)$, then there exist non-negative integers $p, d_j$ and $\mu_j$ such that the unital $C^*$-algebra $C^*(T)$  generated by $T$ is unitarily equivalent to
\[
\oplus_{j=1}^p \bbM_{d_j}(\bbC)^{(\mu_j)}. \]

\bigskip

\begin{prop} \label{prop5.10}
Let $d \ge 2$ and $T \in \cB(\bbC^d)$.    Write $C^*(T) = \oplus_{j=1}^p \bbM_{d_j}(\bbC)^{(\mu_j)}$.  If $p \ge 3$ or if there exists $1 \le j \le p$ such that $d_j \ge 3$, then $\symone(T) \ne \cU(T)$.
\end{prop}

\begin{pf}
Suppose $p \ge 3$.   Then $N := \oplus_{j=1}^p j I_{d_j}^{(\mu_j)} \in C^*(T)$, and $N$ has $p \ge 3$ distinct eigenvalues.  By Proposition~\ref{prop5.09},
\[
\cU(T) \ne \symone(T). \]

\bigskip

If $d_j \ge 3$ for some $1 \le j \le p$, then let $N_i = 0$, $1 \le i \ne j \le p$, and let $N_j := \ttt{diag}(1, 2, \ldots, d_j)$.   Set $N = \oplus_{j=1}^p N_j^{(\mu_j)}$.   Then $N \in C^*(T)$ is normal and $N$ has $d_j \ge 3$ distinct  eigenvalues.  Once again, by Proposition~\ref{prop5.09}, $\symone(T) \ne \cU(T)$.
\end{pf}


\subsection{} \label{sec5.11}
From above, we see that if $2 \le d \in \bbN$, $T \in \cB(\bbC^d)$ and $\cU(T) = \symone(T)$, then $p \le 2$ and $d_j \le 2$, $j = 1, 2$.
We are therefore left with four cases to consider in terms of  the structure of $C^*(T)$.  Let us first deal with the two cases where $\cU(T) \ne \symone(T)$.
\begin{enumerate}
	\item[(i)]
	$p = 1$, $d_1 = 2$.
	
	In this case, $C^*(T) = \bbM_2(\bbC)^{(\mu_1)}$, where $\mu_1 = d/2$, and $T$ is not a normal operator.   Using the canonical shift of basis, we may write an arbitrary element $X \in C^*(T)$ as
	\[
	X = \begin{bmatrix} x_{11} I_{\mu_1} & x_{12} I_{\mu_1} \\ x_{21} I_{\mu_1} & x_{22} I_{\mu_1} \end{bmatrix}, \]
	where $x_{ij} \in \bbC$ is arbitrary.
	We argue by contradiction.  Suppose that $\cU(T) = \symone(T)$.  Then, given $U \in \cU(\bbC^d)$, we can find $J_U \in \symone(\bbC^d)$ such that $U^* T U = J^* T J$, and thus $U^* X U = J^* X J$ for all $X \in C^*(T)$.   In particular, $X (U J^*) = (U J^*) X$ for all $X \in C^*(T)$.
	
	But the commutant of $C^*(T)$ is the set of operators of the form $Y \oplus Y$, where $Y \in \bbM_{\mu_1}(\bbC)$.   Thus
	\[
	V := U J^* = \begin{bmatrix} V_0 & 0 \\ 0 & V_0 \end{bmatrix} \]
	for some unitary operator $V_0 \in \bbM_{\mu_1}(\bbC)$.
	We will now exhibit a specific unitary operator $U$ for which this decomposition is impossible.
	
	Choose $c, s \in (0,1)$ such that $c^2 + s^2 = 1$.   Let $U := \begin{bmatrix} c I_{\mu_1} & - s I_{\mu_1} \\ s I_{\mu_1} & c I_{\mu_1} \end{bmatrix}$.   Then $U$ is unitary.
	
	Suppose that $J = \begin{bmatrix} J_1 & J_2 \\ J_3 & J_4 \end{bmatrix}$ is a symmetry, and that
	\[
	V := U J^* = \begin{bmatrix} V_0 & 0 \\ 0 & V_0 \end{bmatrix}\]
	for some unitary operator $V_0$.
	
	The fact that $J$ is a symmetry implies that $J=J^*$ and thus  $J_3 = J_2^*$ and $J_1, J_2$ are self-adjoint.  Thus
	\[
	\begin{bmatrix} J_1 & J_2 \\ J_3 & J_4 \end{bmatrix} = J = V^* U = \begin{bmatrix} c V_0^* & -s V_0^* \\ s V_0^* & c V_0^* \end{bmatrix}. \]
	That in turn implies that $c V_0^*$ is self-adjoint, and since $0 < c < 1$, that requires that $V_0$ be self-adjoint.

	But since $J_3 = J_2^*$, we also infer that $- s V_0 = (-s V_0^*)^* = (s V_0^*) = s V_0$, a contradiction, as $V_0 \ne 0$.	
	\item[(ii)]
	$p = 2$, $\max(d_1, d_2) = 2$.
	
	Without loss of generality, we may assume that $d_1 = 2$.   Let $N_1 = \begin{bmatrix} 1 & 0 \\ 0 & 2 \end{bmatrix}$ and $N_2 = 3 I_{d_2}$.  Then $N := N_1^{(\mu_1)} \oplus N_2^{(\mu_2)} \in C^*(T)$ is normal and $N$ has three eigenvalues.    By Proposition~\ref{prop5.09}, $\cU(T) \ne \symone(T)$.
\end{enumerate}


\subsection{} \label{sec5.12}
The  two remaining cases are
\begin{enumerate}
	\item[(i)]
	$p=1$, $d_1 = 1$.
	
	In this case, $C^*(T) = \bbC^{(\mu_1)}$, implying that $T = \alpha I_d$ for some $\alpha \in \bbC$.   It is clear that $\cU(T) = \{ T\} = \symone(T)$ in this case.
	\item[(ii)]
	$p =2$, $d_1 = d_2 = 1$.
	
	In this case, there exist $\alpha \ne \beta \in \bbC$ and $1 \le k \le d-1$ such that $T \simeq \alpha I_k \oplus \beta I_{d-k}$.    The proof in this case is somewhat more involved, and will require some auxiliary  definitions and results.
\end{enumerate}	


\begin{defn} \label{def5.13}
Let $1 \le k \le d \in \bbN$ and $T = [t_{ij}] \in \bbM_d(\bbC)$.   If $\Lambda \subseteq \{ 1, 2, \ldots, d\}$ is a set with $k$ elements, then  the $k \times k$ matrix
\[
T_\Lambda := [t_{ij}]_{i, j \in \Lambda} \]
is called the \textbf{principal submatrix} of $T$ of order $k$, determined by $\Lambda$.
\end{defn}


We shall require the following result, which is Theorem~2.20 of~\cite{MarcouxRadjaviZhang2022}.

\begin{thm} \label{thm5.14}
Let $d \in \bbN$ and $T \in \bbM_d(\bbC)$ be a matrix of rank $r$.   Let $\eps > 0$.   Then there exists a unitary matrix $U \in \bbM_d(\bbC)$ satisfying  $\norm U - I_d \norm < \eps$ such that all principal submatrices of $U^* T U$ of order less than or equal to $r$ are invertible.
\end{thm}


\begin{lem} \label{lem5.15}
Let $2 \le d \in \bbN$ and $Y \in \cU(\bbC^d)$.   Suppose that $\bbC^d= \hilb_1 \oplus \hilb_2$, and with respect to this decomposition, let us write $Y = \begin{bmatrix} Y_1  & Y_2 \\ Y_3 & Y_4 \end{bmatrix}$.   Suppose furthermore that $Y_1$ and $- Y_4$ are invertible and positive.
Then $Y \in \symone(\bbC^d)$.
\end{lem}


\begin{pf}
It is not hard to see that by restricting our attention to $\cW^\perp$, where $\cW := \mathrm{span} \{ \ker (Y- I), \ker (Y+I)\}$, it is sufficient to consider the case where $0 < Y_1 < I_{\hilb_1}$ and $0 < -Y_2 < I_{\hilb_2}$, as the restriction of $Y$ to $\cW$ is itself a symmetry.   Furthermore, by interchanging the indices of $\hilb_1$ and $\hilb_2$ and replacing $Y$ by $-Y$ if necessary, we may assume without loss of generality that $\dim\, \hilb_1 \le \dim\, \hilb_2$.

\smallskip

Let us write $M = Y_1$ and note that $0 < M < I_{\hilb_1}$.   Define $N := (I - M^2)^{\frac{1}{2}}$.

The equations $Y^* Y = I = Y Y^*$ imply that
\[M^2+ Y_2 Y_2^* = I_{\hilb_1} = M^2 + Y_3^* Y_3,\]
 or equivalently,
\[
|Y_2^*|^2 = N^2 = |Y_3|^2. \]
It follows from the polar decomposition that there exist co-isometries $W, V: \hilb_2 \to \hilb_1$ such that $Y_2 = N V$ and $Y_3 = W^* N$.  Thus
\begin{align*}
I_{\hilb_2} &= Y_3 Y_3^* + Y_4^2 = W^* N^2 W + Y_4^2, \\
		&= Y_2^* Y_2 + Y_4^2 = V^* N^2 V + Y_4^2.
\end{align*}
In particular,
\[
Y_4 = - W^* M W = - V^* M V. \]
Thus
\[
Y = \begin{bmatrix} M & N V \\ W^* N & - W^* M W \end{bmatrix}. \]
But $Y$ unitary also implies that
\[
0 = Y_3 M + Y_4 Y_2^* = (W^* N) (M) + (-V^* M V)(V^* N). \]
That is,
\[
0 = W^* (NM) - V^* (MN) = (W^* - V^*) MN. \]
Since $MN$ is invertible, $W^* = V^*$, or equivalently, $W = V$.

But then $Y$ is unitary and
\[
Y^* = \begin{bmatrix} M & N W \\ W^*N & -W^* M W \end{bmatrix} = Y, \]
implying that $Y$ is indeed a symmetry.
\end{pf}


\begin{lem} \label{lem5.16}
Let $d \in \bbN$ and suppose that $\bbC^d = \hilb_1 \oplus \hilb_2$, where $\hilb_1, \hilb_2$ are subspaces of $\bbC^d$.   Suppose that $U = \begin{bmatrix} U_1 & U_2 \\ U_3 & U_4 \end{bmatrix} \in \cU(\bbC^d)$.    Then there exist unitary operators $X = X_1 \oplus I_{\hilb_2}$ and $Z = I_{\hilb_1} \oplus Z_4$ in $\cB(\bbC^d)$ such that
\[
X U Z \in \symone(\bbC^d). \]
\end{lem}

\begin{pf}
\noindent{\textsc{Case one:  $U_1$, $U_4$ are invertible.}}\ \ \

Using polar decomposition in the finite-dimensional setting, we may write $U_1 = R_1 |U_1|$ and $U_4 = |U_4^*| R_4$ where $R_1 \in \cU(\hilb_1)$ and $R_4 \in \cU(\hilb_2)$.   Set $X := R_1^* \oplus I_{\hilb_2}$ and $Z := I_{\hilb_1} \oplus - R_4^*$.

Letting $Y_1 = |U_1|$, $Y_2 = -R_1^* U_2 R_4^*, Y_3 = U_3$ and $Y_4 = -|U_4|^*$, we find that
\[
Y := \begin{bmatrix} Y_1 & Y_2 \\ Y_3 & Y_4 \end{bmatrix} = \begin{bmatrix} R_1^* & 0 \\ 0 & I_{\hilb_2} \end{bmatrix} \ \begin{bmatrix} U_1 & U_2 \\ U_3 & U_4 \end{bmatrix} \ \begin{bmatrix} I_{\hilb_1} & 0 \\ 0 & -R_4^* \end{bmatrix} = X U Z. \]

By Lemma~\ref{lem5.15}, $Y \in \symone(\bbC^d)$.

\bigskip

\noindent{\textsc{Case two:  a reduction to the first case.}}\ \ \

By Theorem~\ref{thm5.14}, we can find a sequence $Q_n := \begin{bmatrix} Q_{n1} & Q_{n2} \\ Q_{n3} & Q_{n4} \end{bmatrix}$ of unitary operators such that
\begin{itemize}
	\item{}
	$Q_{n1}$, $Q_{n4}$ are invertible, and
	\item{}
	$U = \lim_n Q_n$.
\end{itemize}
By \textsc{Case one}, for each $n \ge 1$ we can find unitary operators $X_n = X_{n1} \oplus I_{\hilb_2}$ and $Z_n := I_{\hilb_1} \oplus Z_{n4}$ such that
\[
Y_n := X_n U Z_n \in \symone(\bbC^d). \]
By the compactness of the unit ball of $\cB(\bbC^d)$, we can then drop to a subsequence $(Y_{n_k}	)_k$ such that
\[
X := \lim_k X_{n_k}, \ \ \ Z := \lim_k Z_{n_k}\]
both exist.   Clearly $X = X_1 \oplus I_{\hilb_2}$ and $Z = I_{\hilb_1} \oplus Z_4$ for some $X_1 \in \cU(\hilb_1)$ and $Z_4 \in \cU(\hilb_2)$.

But then
\[
Y:= X U Z = \lim_k X_{n_k} Q_{n_k} Z_{n_k} \]
lies in the closure of $\symone(\bbC^d)$.   Since $\symone(\bbC^d)$ is already closed, $Y \in \symone(\bbC^d)$, completing the proof.
\end{pf}


\begin{prop}  \label{prop5.17}
Let $2 \le d \in \bbN$ and suppose that $N \in \cB(\bbC^d)$ is a normal operator with two distinct eigenvalues.   Then $\symone(N) = \cU(N)$.
\end{prop}

\begin{pf}
By hypothesis, we may write $\bbC^d = \hilb_1 \oplus \hilb_2$ such that with respect to this decomposition, $N = \alpha I_{\hilb_1} \oplus \beta I_{\hilb_2}$ with $\alpha \ne \beta \in \bbC$.  

Note that
\[
\uncom(N) = \{ Z_1 \oplus Z_4 :  Z_1 \in \cU(\hilb_1), Z_4 \in \cU(\hilb_2) \}. \]

Let $U \in \cU(\bbC^d)$.   By Lemma~\ref{lem5.16}, there exist $X, Z \in \uncom(N)$ such that $Y := X U Z \in \symone(\bbC^n)$.  By Lemma~\ref{lem5.07}, $\cU(N) = \symone(N)$.
\end{pf}


\bigskip

Together these results imply the following.
\begin{thm} \label{thm5.18}
Let $d \in \bbN$ and $T \in \bbM_d(\bbC)$.   Then $\textsc{Sym}_1(T) = \cU(T)$ if and only if $T$ is normal and has at most two eigenvalues.
\end{thm}

\vskip 1 cm
\vfill\newpage

\begin{center}
\textsc{An ap\'{e}ritif:  $\symtwo$-orbits in low dimensions}
\end{center}


\subsection{} \label{sec5.19}
Having determined the class of operators $T \in \cB(\bbC^d)$ for which $\symone(T)$ coincides with $\cU(T)$, it is reasonable to ask when $\symtwo(T) = \cU(T)$ or $\symthree(T) = \cU(T)$.   While the second of these questions is made particularly difficult because of the fact that we do not have a reasonable characterisation of $\symthree(\bbC^d)$, the fact that $\symtwo(\bbC^d)$ is so tractable suggests that resolving the first of these questions should be much simpler.  Nevertheless, the problem appears to be  difficult, and even for normal operators $N \in \cB(\bbC^d)$, we have only been able to obtain results when $N$ has at most four distinct eigenvalues.



\begin{prop} \label{prop5.21}
Let $T=\oplus_{i=1}^4 \lambda_i I_{d_i}\in \bbM_d(\mathbb{C})$ with $\lambda_1,\lambda_2,\lambda_3,\lambda_4$ being distinct, where $d=\sum_{i=1}^4 d_i$.
Then
\[\mathcal{U}(T)=\textsc{Sym}_2(T).\]

\end{prop}

\begin{pf}
By Lemma~\ref{lem5.07}, it suffices to show that
\[\cU(\bbC^d)=\Gamma_{\cU}(T)\cdot\textsc{Sym}_2(\bbC^d).\]
By Theorem~\ref{thm5.14},
$\cV(\bbC^d)$ is dense in  $\cU(\bbC^d)$, where

\centerline{$\cV(\bbC^d):= \{V: V\in \cU(\bbC^d),$ all principle submatrices of $V$ are invertible$\}$.}
Note that $\Gamma_{\cU}(T)\cdot\textsc{Sym}_2(\bbC^d)\subset \cU(\bbC^d)$, and $\Gamma_{\cU}(T)\cdot\textsc{Sym}_2(\bbC^d)$ is closed.
To finish the proof, we only need to show that
\[\cV(\bbC^d)\subset \Gamma_{\cU}(T)\cdot\textsc{Sym}_2(\bbC^d).\]

Set $\cH_1=\bbC^{d_1}\oplus \bbC^{d_2}$, and
$\cH_2=\bbC^{d_3}\oplus \bbC^{d_4}$.
Set $T_1=\lambda_1 I_{d_1}\oplus \lambda_2 I_{d_2}$, and $T_2=\lambda_3 I_{d_3}\oplus \lambda_4 I_{d_4}$.
Then $T=T_1\oplus T_2$.

Given $V\in \cV(\bbC^d)$, with respect to $\bbC^d=\cH_1\oplus \cH_2$,
we write
\[V=\begin{bmatrix} V_1& V_2 \\ V_3& V_4 \\ \end{bmatrix},\]
where $V_1$ and $V_3$ are invertible.  Using polar decomposition in finite-dimensional setting, we may write $V_1=R_1|V_1|$, $V_4=R_4|V_4|$,
where $R_1\in \cU(\cH_1)$, $R_4\in \cU(\cH_2)$.
Now set
\[Y:=\begin{bmatrix} Y_1& Y_2 \\ Y_3& Y_4 \\ \end{bmatrix}=\begin{bmatrix} R_1^*& 0 \\ 0& -R_4^* \\ \end{bmatrix}\begin{bmatrix} V_1& V_2 \\ V_3& V_4 \\ \end{bmatrix}=\begin{bmatrix} |V_1|& R_1^*V_2 \\ -R_4^*V_3& -|V_4| \\ \end{bmatrix}.\]
By Lemma~\ref{lem5.15}, $Y\in \textsc{Sym}_1(\bbC^d)$.

By Theorem~\ref{thm5.18},
\[\textsc{Sym}_1(T_i)=\cU(T_i),~~i=1,2.\]
Then by Lemma~\ref{lem5.07}, $\cU(\cH_i)=\Gamma_{\cU}(T_i)\cdot \textsc{Sym}_1(\cH_i)$, $i=1,2$.
Hence, we could choose $W_i\in \Gamma_{\cU}(T_i), J_i\in \textsc{Sym}_1(\cH_i)$, $i=1,2$,  such that $R_1=W_1J_1$ and $-R_4=W_2J_2$.
Observe that $W_1\oplus W_2\in \Gamma_{\cU}(T)$ and $J_1\oplus J_2\in \textsc{Sym}_1(\bbC^d)$.
Now
\[V=\begin{bmatrix} R_1& 0 \\ 0& -R_4 \\ \end{bmatrix}Y=\begin{bmatrix} W_1& 0 \\ 0& W_2 \\ \end{bmatrix}\begin{bmatrix} J_1& 0 \\ 0& J_2 \\ \end{bmatrix}Y\in \Gamma_{\cU}(T)\cdot\textsc{Sym}_2(\bbC^d).\]
\end{pf}


Using a similar argument, we can obtain the following corollary.

\begin{cor} \label{cor5.22}
Let $T=\oplus_{i=1}^3 \lambda_i I_{d_i}\in \cB(\mathbb{C}^d)$ with $\lambda_1,\lambda_2,\lambda_3$ being distinct, where $d=\sum_{i=1}^3 d_i$.
Then
\[\mathcal{U}(T)=\textsc{Sym}_2(T).\]
\end{cor}


\subsection{} \label{sec5.33}
Despite considerable effort, we have not yet resolved the question of whether or not $\symtwo(N) = \cU(N)$ (or even $\symthree(N) = \cU(N)$) when $N$ is a normal (and thus diagonalisable) operator in $\cB(\bbC^d)$ when $d \ge 5$.  Indeed, suppose that $D$ is a diagonal operator with respect to some orthonormal basis $\cE := \{ e_k: 1 \le k \le d\}$, and to simplify the problem even further, consider the case where $D$ admits $d$ distinct eigenvalues.   In this case, the unitary commutant $\Gamma_\cU(D) = \cD(\bbC^d)$, the set of diagonal operators with respect to $\cE$.   By Lemma~\ref{lem5.07}, the question of whether $\symtwo(D) = \cU(D)$ is equivalent to the question:  given $U \in \cU(\bbC^d)$, can we find a diagonal, unitary operator $B$ such that the non-real eigenvalues of $B U$ appear in conjugate pairs?

To see why this seemingly innocuous question may be more difficult than it seems, we are tempted to draw a parallel with the well-known conjecture of Horn~\cite{Horn1962} regarding the spectrum of the sum of two hermitian operators, which was eventually resolved by Klyachko~\cite{Klyachko1998} and independently by Knuston and Tao~\cite{KnutsonTao1999}.     Given unitary operators $U$ and $V$ in $\cU(\bbC^d)$, we may find hermitian matrices $H, K \in \bbM_d(\bbC)$ such that $U = e^{i H}$ and $V = e^{i K}$.   Of course, the equation $U V = e^{i (H+K)}$ does not hold unless $U$ and $V$ commute (in which case the answer to our original problem is easily seen to be ``yes"), but the analogy suggests that calculating the possible spectra of operators of the form $U V$ in general may be akin to calculating the possible spectra of $H + K$.  We hope that the solution to the multiplicative version of this problem for unitaries will not require thirty-six years to solve, as did the Horn Conjecture.

%
%
%
%


\bibliographystyle{plain}


\end{document}